\documentclass[preprint, 12pt]{elsarticle}
\usepackage{amsfonts}
\usepackage{amsmath}

\newcommand\R{{\mathbb{R}}}

\pdfoutput=1

\usepackage{tikz}

\usepackage{natbib}
\usepackage{bibentry}

\usepackage{subcaption}

\begin{document}

\newtheorem{theorem}{Theorem}[section]
\newtheorem{lemma}[theorem]{Lemma}
\newtheorem{Pro}[theorem]{Proposition}
\newtheorem{definition}[theorem]{Definition}
\newtheorem{proposition}[theorem]{Proposition}
\newtheorem{remark}[theorem]{Remark}
\newtheorem{corollary}[theorem]{Corollary}

\title{Convergence of spectra of uniformly fattened open book structures}

\author{James E. Corbin}
\address{Department of Mathematics\\
Texas A\& M University\\
College Station, TX}
\ead{jamesedwcorbin@gmail.com}



\begin{abstract}
We consider a compact $C^\infty$-stratified $2D$ variety $M$ in $\R^3$ 
and its $\epsilon$--neighborhood $M_\epsilon$, which we call a ``fattened open book structure.'' 
Assuming absence of zero-dimensional strata, i.e. ``corners'', 
we show that the (discrete) spectrum of the Neumann Laplacian in $M_\epsilon$ converges when $\epsilon\to 0$ to the spectrum of
a differential operator on $M$.

Similar results have been obtained before for the case of fattened graphs, i.e. $M$ being one-dimensional. 
In the case of a $2D$ smooth submanifold $M$, the problem has been studied well. 
However, having singularities along strata of lower dimensions significantly complicates considerations. 
As in the quantum graph case, such considerations are triggered by various applications.
\end{abstract}


\date{\today}
\maketitle

\section{Introduction}

\begin{figure}[ht!]
\begin{center}
\includegraphics{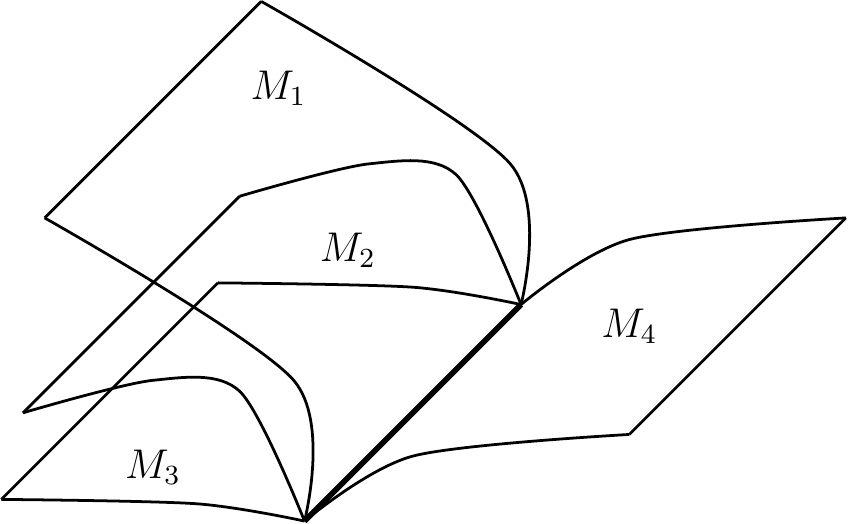}
\end{center}
		 \caption{An open book structure with ``\textbf{pages}'' $M_k$ meeting at a ``\textbf{binding}.''}
		 \label{F:OBS}
\end{figure}

This article is the fulfillment of the announced results in \cite{CoKu}.
We consider a compact $C^\infty$-stratified
\footnote{We do not provide here the general definition of what is called Whitney stratification, 
		see e.g. \cite{Arnold,Lu,DifTop,Whit,strat}, resorting to a simple description through local models.} 
$2D$ sub-variety $M$ in $\R^3$ without zero-dimensional strata, i.e.
$M$ locally (in a neighborhood of any point) looks like either a smooth submanifold 
or like an ``open book'' with smooth two-dimensional ``pages'' meeting transversely 
along a common smooth one-dimensional ``binding,'' see Fig. \ref{F:OBS}.
The $2D$-strata need neither be contractible nor orientable.

Clearly, any compact smooth submanifold of $\R^3$ (with or without a boundary) qualifies as 
an open book structure with a single page. 
Another example of such structure is shown in Fig. \ref{F:spheres}:
\begin{figure}[ht!]
		\begin{center}
\includegraphics{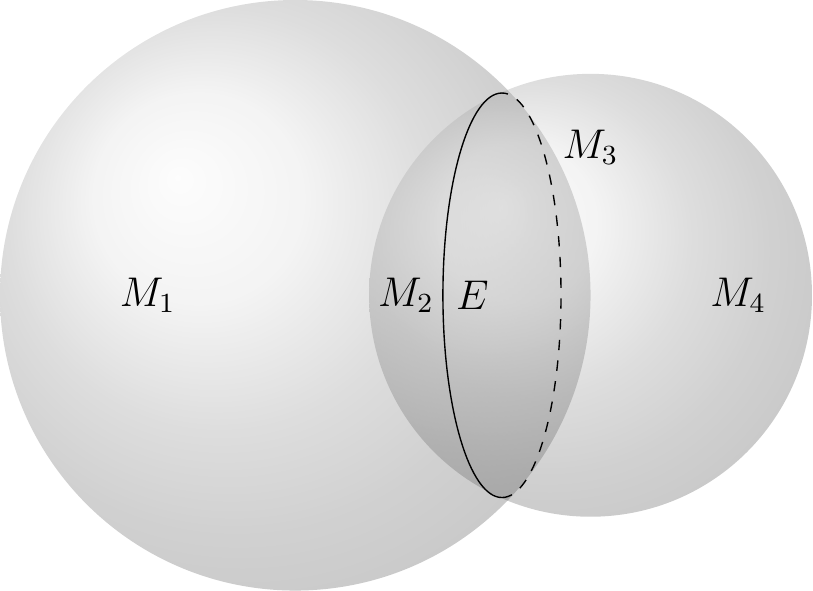}
		\end{center}
		\caption{A transverse intersection of two spheres yields an open book structure with
		four pages and a circular binding. 
		The requirement of absence of zero-dimensional strata prohibits
		adding a third sphere with a generic triple intersection. 
		Tangential contacts of spheres are also disallowed.}
		\label{F:spheres}
\end{figure}

We then consider a ``fattened'' version $M_\epsilon$ of $M$, which is
an (appropriately defined) $\epsilon$--neighborhood of $M$, which we call a ``fattened open book structure.''

Consider now the Laplace operator $-\Delta$ on the domain $M_\epsilon$ with Neumann boundary conditions (the ``\textbf{Neumann Laplacian}'').
We denote this operator\footnote{ Throughout this paper $\epsilon$-dependent spaces, functions, 
		and coefficients will carry an $\epsilon$ subscript or superscript.} $A_\epsilon$.
As a (non-negative) elliptic operator on a compact manifold, 
it has discrete finite multiplicity spectrum $\lambda^\epsilon_n:=\lambda_n(A^\epsilon)$ with the only accumulation point at infinity.
The main result of this work is that when $\epsilon\to 0$, each eigenvalue $\lambda_n^\epsilon$ 
converges to the corresponding eigenvalue $\lambda_n$ of an operator $A$ on $M$, 
which acts as $-\Delta_M$ (2D Laplace-Beltrami) on each 2D stratum (\textbf{page}) of $M$, 
with appropriate junction conditions along 1D strata (\textbf{bindings}).

Similar results have been obtained previously for the case of fattened graphs 
(see \cite{KZ,RS2, Fr, FW}, as well as books \cite{BerKuc,Post} and references therein), 
i.e. $M$ being one-dimensional with some weak results for open books \cite{FWopb}. 
They have been triggered by various applications \cite{DellAn06, ES, FK1,WRM, RS2, RS3, RS5, RSS}. 
In the case of a \emph{smooth} submanifold $M\subset \R^3$, 
the problem is not that hard and has been studied well under a variety of ``hard'' and ``soft'' 
constraints set near $M$ (see, e.g. \cite{FH,Grieser,WRM}). 
However, having singularities along strata of lower dimensions significantly complicates considerations, 
even in the quantum graph case \cite{Grieser, DellAn06, DellAn07, DellAn15, KuISAAC, KZ, KZ2, WRM, RS2, S, EH}.

The paper is structured as follows: Section \ref{S:notions} contains the descriptions of the main objects: open book structures and their fattened versions, 
the limit operator $A$, etc. 
The next Section \ref{S:result} contains formulation of the main result. 
The proofs are provided in Section \ref{S:proof_I}, and
Section \ref{S:remarks} contains the final remarks and discussions.

In this article the results are obtained under two restrictions: that the width of the fattened domain shrinks ``with the same speed'' around all strata and the binding of the book does not have a cusp.
The more complex case of slower shrinkage of the neighborhoods of lower dimensional strata, which leads to phase transition phenomena (see \cite{KZ2,Post} for the quantum graph case), will be considered elsewhere. 
The same applies to the even more complex case of the presence of corners.

\section{The main notions}\label{S:notions}
Here we introduce the main geometric objects and to be studied.
\subsection{Open book structures and their fattened counterparts}\label{SS:OBS}

For the purpose of this paper, we deal with a particular family of stratified varieties called open book structures.
\footnote{One can find open book structures in a somewhat more general setting being discussed in algebraic topology literature, e.g. in \cite{Ranic,Winkel}.}

\begin{definition} \label{D:M}
We call $M$, a connected compact $C^\infty$ stratified two-dimensional variety in $\mathbb{R}^3$, an \textbf{open book structure} if:
\begin{itemize}
\item zero dimensional strata are absent;
\item $M$ is composed of finitely many smooth $2D$ strata $\{M_k\}$ ($k \leq n_M$) (open smooth surfaces) called \textbf{pages} and smooth $1D$ strata $\{E_m\}$ ($m \leq n_E$) (close smooth curves) called \textbf{bindings}.
\item The pages are transverse at the bindings.
\end{itemize}
\end{definition}

Simply put, to say $M$ is a stratified surface means that it consists of finitely many connected, 
compact smooth submanifolds (with or without boundary) of $\R^3$, called \textbf{strata}, 
of dimensions two, one, and zero  
such that they may only intersect along their boundaries and each stratum's boundary is the union of some lower dimensional strata \cite{strat}.
Point strata, zero dimensional strata, are omitted from consideration in this article.
Additionally, We assume that the strata intersect at their boundaries transversely.


%
%


Next, we build a new domain by ``fattening'' $M$ by considering its $\epsilon$-neighborhood $M_\epsilon$.
In the planned sequel to this article we consider a more general variable width neighborhood, 
with a smooth function $w(x)$ on $M$ controlling the variable width\footnote{When $w=1$, this boils down to the standard neighborhoods $M_\epsilon$.}.
We denote the ball of radius $r$ at $x$ as $B(x, r)$.

\begin{definition} \label{D:Me}
		The \textbf{uniformly fattened domain} $M_\epsilon$ ($\epsilon>0$) over an open book structure $M$ consists of all points of distance $\epsilon$ from $M$.
		I.e.
		\begin{equation}
				M_\epsilon = \bigcup_{x \in M} B(x, \epsilon).
		\end{equation}
		The similar notation $R_\epsilon$ will be used for the fattened version of any subset $R\in\R^3$.
\end{definition}

The following statement is rather obvious:
\begin{lemma} \label{L:epsilon0}
	There exists $\epsilon>0$ so small that for any two points $x_1, x_2\in M$ outside of an 
	$\epsilon_0$-neighborhood of the bindings, the closed intervals of radius $\epsilon_0$ normal to $M$ at these points do not intersect.
\end{lemma}

This ensures that for $\epsilon\leq \epsilon_0$, the $\epsilon$-fattened neighborhoods do not form a connecting bridge between two points that are otherwise far away from each other along $M$.
Henceforth we assume $\epsilon \leq \epsilon_0$, and our goal is to analyze the behavior of operators on $M_\epsilon$ in the $\epsilon \to 0$ limit.
Therefore notations like $O(\epsilon)$ or $o(1)$ are understood with respect to that limit.


\subsection{Local structure} \label{SS:local structure}

\begin{figure}[ht!]
	\begin{center}
		\begin{subfigure} {0.8\textwidth}
\includegraphics{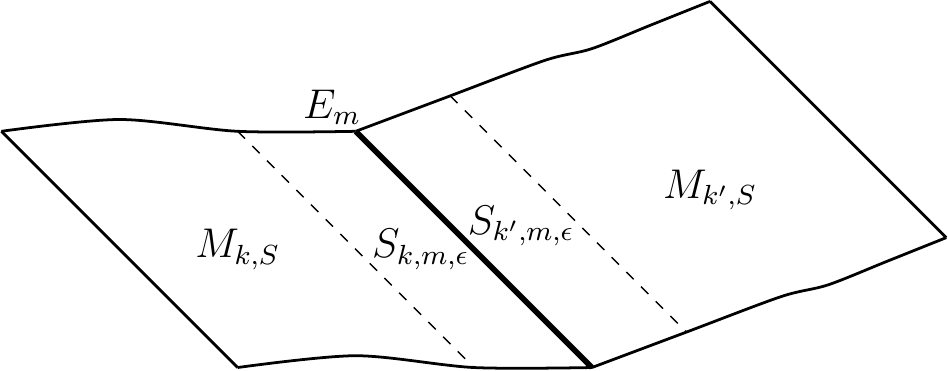}
		\end{subfigure}
		\begin{subfigure}{0.7\textwidth}
\includegraphics{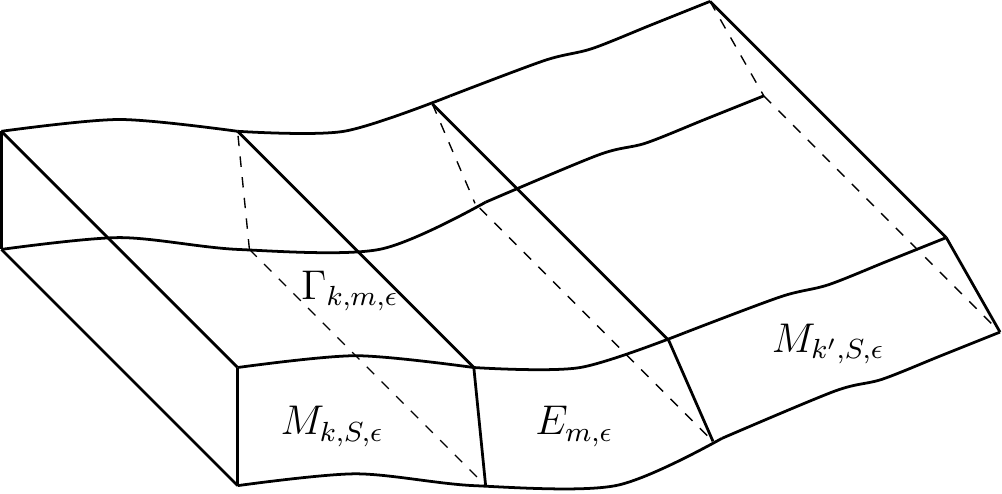}
		\end{subfigure}
\caption{A neighborhood of a binding and the corresponding uniformly fattened neighborhood. }
\label{F:MandMe}
	\end{center}
\end{figure}

For any binding $E_m$, the parts of the adjacent pages $M_k$ that are $O(\epsilon)$--close to $E_{m}$ will be called \textbf{sleeves} and denoted $S_{k,m,\epsilon}$. 
More precisely,

\begin{definition} \label{D:Ske}
	Let $M$ be an open book structure with $n_E$ bindings. 
	Let $\{a_m\}_{m\leq n_E}$ denote a finite set of positive numbers independent of $\epsilon$.
	The sleeve $S_{k,m,\epsilon}$ on page $M_k$ at $E_m$ is defined as
	\begin{equation}
		S_{k,m,\epsilon} := \{ x \in M_k : \mathsf{dist}_{M_k}(x, E_m) < a_m \epsilon \},
	\end{equation}
	where $\mathsf{dist}_{M_k}(x, E_m)$ denotes the geodesic distance from $E_m$ to $x$ on $M_k$ (see Fig. \ref{F:MandMe}).

	We will use the following shorthand notation for the page without its sleeves:
	$$M_{k,S}:= M_k \backslash \bigcup_m S_{k,m,\epsilon}.$$
\end{definition}

The next statement is easy to establish due to the non-transverse nature of pages' intersections:
\begin{lemma} \label{L:Ske}
	Under appropriate choice (which we will fix) of $\{a_m\}$, the $\epsilon$-neighborhoods of $M_{k,S}$ does not intersect each other for different values of $k$ and any binding $E_m$.
\end{lemma}

\begin{definition}\label{D:normal}
	Assuming a choice of orientation of $M_k$, we denote the unit normal vector to $M_k$ at a point $x\in M_k$ as $\mathcal{N}_k(x)$.
	If $M_k$ is non-orientable, a local choice of normal orientation will be sufficient for our purposes.
	We denote by $\mathcal{I}_{\mathcal{N}_k(x),\epsilon}$ the interval of the normal to $M_k$ at $x$ consisting of points at distance less than $\epsilon$ from $y$.
\end{definition}

The fattened page $M_{k, S, \epsilon}$ is thus foliated into normal fibers $\mathcal{I}_{\mathcal{N}_k(x),\epsilon}$.
The latter foliation will be used to define the local averaging operator on $M_{k, S, \epsilon}$ in Section \ref{SS:ir_I}.

\begin{definition} \label{D:Ee}
	The fattened binding $E_{m,\epsilon}$ about $E_m$ is the union of the $\epsilon$-neighborhood of $E_m$ 
	and the $2 \epsilon$ width normal fibers over the sleeves $S_{k,m,\epsilon}$ (see Fig. \ref{F:MandMe}):
	\begin{equation}
		E_{m,\epsilon} := \bigcup_{x \in E_m} B(x, \epsilon) \bigcup \big( \bigcup_{k; x \in S_{k,m,\epsilon}} \mathcal{I}_{\mathcal{N}_k( x ), \epsilon} \big)
	\end{equation}

	We can also define a \textbf{cross-section} $\omega_{m, \epsilon}(x)$.
	For a point $x$ in $E_m$, $N_{x}$ is the normal plane of $E_m$ at $x$, an affine subspace of $\mathbb{R}^3$.
	The cross-section $\omega_{m, \epsilon}(x)$ is the connected component of the intersection of $N_x$ with $M_{\epsilon} \backslash \bigcup_k M_{k, S, \epsilon}$ containing $x$.

	The fattened binding can also be defined as the union of these cross-sections:
	\begin{equation}
		E_{m, \epsilon} := \bigcup_{x \in E_m} \omega_{m, \epsilon}(x).
	\end{equation}
\end{definition}

\subsection{Quadratic forms and operators}\label{SS:qforms}

We adopt the standard notation for Sobolev spaces (see, e.g. \cite{Ma}). 
Thus, $H^1(\Omega)$ denotes the space of square integrable with respect to the 
Lebesgue measure functions on a domain $\Omega \subset \mathbb{R}^n$ with square integrable first order weak derivatives.

\begin{definition} \label{D:Qe}
		Let $Q_\epsilon$ denote the closed non-negative quadratic form with domain $H^1(M_\epsilon)$, given by
		\begin{equation}
			 Q_{\epsilon} (u ) = \int_{M_{\epsilon}} | \nabla u |^2 \, dM_\epsilon
		\end{equation}
		We also refer to $Q_{\epsilon}(u)$ as the \textbf{energy} of $u$.
\end{definition}

This form is associated with a unique self-adjoint operator $A_\epsilon$ in $L_2(M_\epsilon)$. 
The following statement is standard (see, e.g. \cite{EE,Ma}):

\begin{Pro} \label{P:QeAe}
		The form $Q_\epsilon$ generates the \textbf{Neumann Laplacian} on $M_\epsilon$. 
		I.e. $A_\epsilon = -\Delta$ with its domain consisting of functions in  $H^2(M_\epsilon)$ whose normal derivatives at the boundary $\partial M_\epsilon$ vanish.

		Its spectrum $\sigma(A_\epsilon)$ is discrete and non-negative.
\end{Pro}

We equip $M$ with the surface measure $dM$ induced from $\R^3$.
Proceeding to the open book structure, we define the energy on $M$ as follows:

\begin{definition} \label{D:Q}
		Let $Q$ be the closed, non-negative quadratic form (\textbf{energy}) on $L_2(M)$ given by
		\begin{equation}
			Q(u) = \sum_k \int_{M_k} | \nabla_{M_k} u |^2 \, dM
		\end{equation}
		with domain $\mathcal{G}^1$ consisting of functions $u$ for whose $Q(u)$ is finite and that are continuous across the bindings between pages $M_k$ and $M_{k'}$:
		\begin{equation}\label{E:cont}
			u|_{\partial M_k \cap E_m} = u|_{\partial M_{k'} \cap E_m}.
		\end{equation}
		Here $\nabla_{M_k}$ is the gradient along $M_k$ and restrictions in (\ref{E:cont}) to the binding $E_m$ coincide as elements of $H^{1/2}(E_m)$.
\end{definition}

Unlike the fattened graph case, by the Sobolev embedding theorem \cite{EE} 
the restriction to the binding is not continuous as an operator from $\mathcal{G}^1$ to $C(E_m)$, 
it only maps to $H^{1/2}(E_m)$. 
This distinction significantly complicates the analysis of fattened stratified surfaces in comparison with fattened graphs.

\begin{Pro} \label{Ae}
	The operator $A$ associated with the quadratic form $Q$ acts on each $M_k$ as
		\begin{equation}
				A u :=
				-\Delta_{M_k} u,
		\end{equation}
		with the domain $\mathcal{G}^2$ consisting of functions on $M$ such that
		the following conditions are satisfied:
		\begin{equation}
				||u||_{L_2(M)}^2 +  ||A  u||_{L_2(M)}^2  < \infty,
		\end{equation}
		continuity across common bindings $E_m$ of pairs of pages $M_k,M_{k'}$:
		\begin{equation}
				u|_{\partial M_k \cap E_m} = u|_{\partial M_{k'} \cap E_m},
		\end{equation}
		\textbf{Kirchhoff condition} at the bindings:
		\begin{equation} 
				\sum_{k: \partial M_k \supset E_m  }  D_{\nu_k} u (E_m) = 0.
		\end{equation}
		Here $-\Delta_{M_k}$ is the Laplace-Beltrami operator on $M_k$ and $D_{\nu_k}$ denotes the normal derivative to $\partial M_k$ along $M_k$.

		The spectrum of $A$ is discrete and non-negative.
\end{Pro}

The proof is simple, standard, and similar to the graph case. We thus omit it.

\begin{definition}
		For a real number $\Lambda$ not in the spectrum of $A_\epsilon$, 
		we denote by $\mathcal{P}_{\Lambda}^\epsilon$ the spectral projector 
		of $A_{\epsilon}$ in $L_2(M_\epsilon)$ onto the spectral subspace 
		corresponding to the half-line $\{\lambda\in\R\, |\, \lambda < \Lambda\}$.

		Similarly, $\mathcal{P}_{\Lambda}$ denotes the analogous spectral projector for $A$.
		We then denote the corresponding (finite dimensional) spectral subspaces 
		as $\mathcal{P}_{\Lambda}^\epsilon L_2(M_\epsilon)$ 
		and $\mathcal{P}_{\Lambda} L_2(M)$ for $M_\epsilon$ and $M$ respectively.
\end{definition}

\begin{Pro} \label{P:proj}
	Functions from these (finite-dimensional) spectral subspaces satisfy the ``reverse" embedding inequality. 
	Namely, if $u \in \mathcal{P}_{\Lambda}^\epsilon L_2(M_\epsilon)$ 
	and $\Lambda \notin \sigma(A_\epsilon)$ then $u \in H^1(M_\epsilon)$ with
		\begin{equation}
			||u ||_{H^1(M_\epsilon)}^2 \leq (1+\Lambda) ||u||^2_{L_2(M_\epsilon)}
		\end{equation}
		and similarly $u \in \mathcal{P}_{\Lambda}  L_2(M )$ and $\Lambda \notin \sigma(A)$
		\begin{equation}
			||u ||_{\mathcal{G}^1}^2 \leq (1+\Lambda) ||u||^2_{L_2(M )}
		.
		\end{equation}
\end{Pro}

\textbf{Proof:}
The proof is a simple application of spectral theory, and so it is omitted.

\section{Formulation of the main result}\label{S:result}

We denote the non-decreasingly ordered eigenvalues of $A$ as $\{\lambda_n\}_{n\in \mathbb{N}}$ 
and those of $A_\epsilon$ as $\{\lambda_n^\epsilon \}_{n \in \mathbb{N}}$.  

\begin{definition}
		We say the operators $A_\epsilon$ \textbf{converge in spectra to} $A$ as $\epsilon$ tends to zero if for each $n$
		$$ | \lambda_n - \lambda_n^\epsilon | = o(1), $$

		where $o(1)$ \underline{is not necessarily uniform with respect to $n$}.
\end{definition}

We now introduce two families of operators needed for the formulation and proof of the main result.

\begin{definition} \label{D:Je}
		A family of linear operators $J_{\epsilon}$
		from $H^1(M_{\epsilon})$ to $\mathcal{G}^1$ is called \textbf{averaging operators} 
		if for any $\Lambda \notin \sigma(A_\epsilon)$ there is an $\epsilon_0$ such that for all $\epsilon \in (0, \epsilon_0]$ the following conditions are satisfied:
		\begin{itemize}
				\item For $u \in \mathcal{P}_{\Lambda}^\epsilon L_2(M_{\epsilon})$, $J_\epsilon$ is ``nearly an isometry'' from $L_2(M_\epsilon)$ to $L_2(M)$ with an $o(1)$ error, i.e.
				\begin{equation}\label{E:J1}
				\bigg| \, \big| \big| u \big| \big|_{L_2(M_{\epsilon})}^2 - \big| \big| J_\epsilon u \big| \big|_{L_2(M)}^2 \, \bigg| \leq o(1) \big| \big| u \big| \big|_{H^1(M_{\epsilon})}^2
				\end{equation}
				where $o(1)$ is uniform with respect to $u$.
				\item For $u \in \mathcal{P}_{\Lambda}^\epsilon L_2(M_{\epsilon})$, $J_\epsilon$ asymptotically ``does not increase the energy,'' i.e.
				\begin{equation}\label{E:J2}
				 Q( J_\epsilon u) - Q_{\epsilon} ( u )  \leq o(1) Q_{\epsilon}(u)
				\end{equation}
				where $o(1)$ is uniform with respect to $u$.
		\end{itemize}
\end{definition}

\begin{definition} \label{D:Ke}
		A family of linear operators $K_{\epsilon}$ from 
		$\mathcal{G}^1$ to $H^1(M_{\epsilon})$ is called \textbf{extension operators} 
		if for any $\Lambda \notin \sigma(A)$ there is an $\epsilon_0$ such that 
		for all $\epsilon \in (0, \epsilon_0]$ the following conditions are satisfied:
		\begin{itemize}
				\item  For $u \in  \mathcal{P}_{\Lambda} L_2(M)$, $K_\epsilon$ is ``nearly an isometry'' from $L_2(M)$ to $L_2(M_\epsilon)$ with $o(1)$ error, i.e.
				\begin{equation}\label{E:K1}
				\bigg| \, \big| \big| u \big| \big|_{L_2(M)}^2 - \big| \big| K_\epsilon u \big| \big|_{L_2(M_{\epsilon})}^2 \, \bigg| \leq o(1) \big| \big| u \big| \big|_{\mathcal{G}^1}^2
				\end{equation}
				where $o(1)$ is uniform with respect to $u$.
				\item For $u \in \mathcal{P}_{\Lambda} L_2(M) $, $K_\epsilon$ asymptotically ``does not increase'' the energy, i.e.
				\begin{equation}\label{E:K2}
				Q_{\epsilon}( K_\epsilon u)- Q ( u )   \leq o(1) Q(u)
				\end{equation}
				where $o(1)$ is uniform with respect to $u$.
		\end{itemize}
\end{definition}

Existence of such averaging and extension operators is known to be 
sufficient for spectral convergence of $A_\epsilon$ to $A$ (see \cite{Post}). 
We now precisely formulate this in our situation.

\begin{theorem} \label{T:spec}
		Let $M$ be an open book structure as in Definition \ref{D:M} and its fattened partner $\{M_\epsilon\}_{\epsilon \in (0, \epsilon_0]}$ as in Definition \ref{D:Me}.
		Let $A$ and $A_\epsilon$ be operators on $M$ and $M_\epsilon$ as in Definitions \ref{P:QeAe} and \ref{Ae}.

		Suppose there exist averaging operators $\{J_\epsilon\}_{\epsilon \in (0, \epsilon_0]}$ 
		and extension operators $\{K_\epsilon\}_{\epsilon \in (0, \epsilon_0]}$ as stated in Definitions \ref{D:Je} and \ref{D:Ke}.

		Then, for any $n$ $$\lambda_n(A_\epsilon) \mathop{\rightarrow}_{\epsilon \to 0} \lambda_n(A).$$
\end{theorem}

We start with the following standard (see, e.g. \cite{ReSi}) min-max characterization of the spectrum.
\begin{proposition}\label{P:Rayleigh}
Let $B$ be a self-adjoint non-negative operator with discrete spectrum of finite multiplicity and $\lambda_n(B)$ be its eigenvalues listed in non-decreasing order. Let also $q$ be its quadratic form with the domain $D$. Then
\begin{equation}\label{E:minmax}
\lambda_n(B)=\mathop{min}\limits_{W\subset D}\,\,\mathop{max}\limits_{x\in W\setminus \{0\}}\frac{q(x,x)}{(x,x)},
\end{equation}
where the minimum is taken over all $n$-dimensional subspaces $W$ in the quadratic form domain $D$
\end{proposition}

\textbf{Proof of Theorem \ref{T:spec}} now employs Proposition \ref{P:Rayleigh} and the averaging and extension operators $J,\,K$ to ``transplant'' the test spaces $W$ in  (\ref{E:minmax}) between the domains of the quadratic forms $Q$ and $Q_\epsilon$.

Let us first notice that due to the definition of these operators (the near-isometry property), for any fixed finite-dimensional space $W$ in the corresponding quadratic form domain, for sufficiently small $\epsilon$ the operators are injective on $W$ and thus preserve its dimension. Since we are only interested in the limit $\epsilon\to 0$, we will assume below that $\epsilon$ is sufficiently small for these operators to preserve the dimension of $W$. Thus, taking also into account the inequalities (\ref{E:J1})-(\ref{E:K2}), one concludes that on any fixed finite dimensional subspace $W$ one has the following estimates of Rayleigh ratios:
\begin{equation}\label{E:RrJ}
\dfrac{ Q(J_\epsilon u) }{ || J_\epsilon u ||^2_{L_2(M)} } \leq \big(1 + o(1) \big)  \dfrac{ Q_\epsilon( u ) }{ || u||^2_{L_2(M_\epsilon)}}
\end{equation}

\begin{equation}\label{E:RrK}
\dfrac{ Q_{\epsilon} ( K_\epsilon u) }{ || K_\epsilon u ||_{L_2(M_\epsilon)}^2} \leq \big( 1 + o(1) \big)  \dfrac{ Q(u) }{ || u ||_{L_2(M)}^2}
\end{equation}

Let now $W_n\subset \mathcal{G}^1$ and $W_n^\epsilon\subset H^1(M_\epsilon)$ be $n$, such that
\begin{equation}\label{E:minmax}
\lambda_n=\mathop{max}\limits_{x\in W_n\setminus \{0\}}\frac{Q(x,x)}{(x,x)},
\end{equation}
and
\begin{equation}\label{E:minmax}
\lambda^\epsilon_n=\mathop{max}\limits_{x\in W^\epsilon_n\setminus \{0\}}\frac{Q_\epsilon(x,x)}{(x,x)},
\end{equation}
Due to the min-max description and inequalities (\ref{E:RrJ}) and (\ref{E:RrK}), one gets

\begin{equation}
\lambda_n \leq
\sup_{u \in J_\epsilon(W_n^\epsilon)} \dfrac{ Q(J_\epsilon u) }{ || J_\epsilon u ||^2_{L_2(M)}}
 \leq \big( 1 + o(1) \big)  \lambda_n^\epsilon,
\end{equation}
and
\begin{equation}
\lambda_n^\epsilon \leq
\sup_{u \in K_\epsilon(W_n)} \dfrac{ Q_\epsilon(K_\epsilon u) }{ || K_\epsilon u ||^2_{L_2(M_\epsilon)}}
\leq \big( 1 + o(1) \big)  \lambda_n
\end{equation}
Thus, $ \lambda_n  - \lambda_n^\epsilon = o(1)$, which proves the theorem. $\Box$

We will construct the required averaging and extension operators, which then will lead to the main result of this text:

\begin{theorem} \label{T:main}
		Let $M$ be an open book structure as in Definition \ref{D:M} and its fattened partner $\{M_\epsilon\}_{\epsilon \in (0, \epsilon_0]}$ as in Definition \ref{D:Me}.
		Let $A$ and $A_\epsilon$ be operators on $M$ and $M_\epsilon$ as in Definitions \ref{P:QeAe} and \ref{Ae}. 
		There exist averaging operators $\{J_\epsilon\}_{\epsilon \in (0, \epsilon_0]}$ 
		and extension operators $\{K_\epsilon\}_{\epsilon \in (0, \epsilon_0]}$ as stated in Definitions \ref{D:Je} and \ref{D:Ke}.

		Thus, for any $n$ $$\lambda_n(A_\epsilon) \mathop{\rightarrow}_{\epsilon \to 0} \lambda_n(A).$$
\end{theorem}


\section{The Proof of the Main Result (Theorem \ref{T:main})}\label{S:proof_I}

In order to define these averaging and extension operators, we must first consider the different local geometries of $M$.
We define a local averaging operator on each of the fattened strata and a local extension operator from each of the pages into $M_\epsilon$.
Then we find a way to reconcile these local operators defined on different geometries.
 This is somewhat similar to the analysis on the fattened graph; however, different embedding theorems in dimensions higher than $1$ require a more careful analysis than in the graph case.

\subsection{Fattened Binding Geometry} \label{SS:FBG}

In this subsection we describe the geometry of the fattened binding and, in particular, specify the lengths $a_m$.
We describe carefully the geometry in order for the domain to admit a suitable partition of unity.
This partition of unity is chosen as to allow good estimates with regards to $\epsilon$ dependence on the norms of trace and extension operators.

\begin{figure}
\centering
\includegraphics{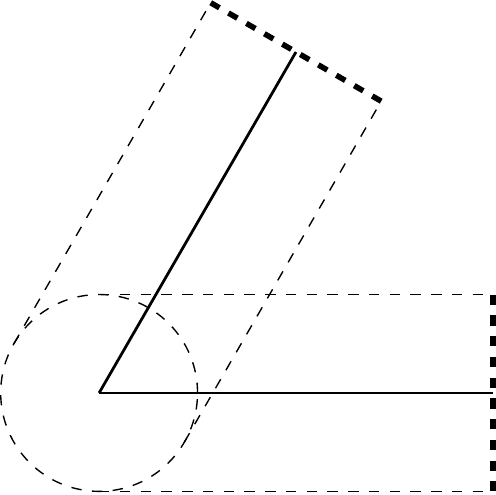}
\caption{A cross-section of a uniformly fattened binding neighborhood. Dashed lines denote the boundary of a fattened stratum. Thickest dashed lines denote the cross-section of the boundary $\Gamma_{ k, m, \epsilon}$ between the fattened binding and fattened page.}
\label{F:Eme_I}
\end{figure}

\begin{definition} 
\label{D:a_I}
Let $M$ be an open book structure.
Let $\theta_{m,k,k'}(x)$ be the (smaller) angle between two tangent vectors normal to two intersecting page boundaries $\partial M_k$ and $\partial M_{k'}$ at $x \in E_m$.
The sleeve width $a_{m }$ ($m \leq n_E$) (see Fig. \ref{F:Eme_I}) is
\begin{equation}
a_{m } = \begin{cases}
 \max_{x \in E_m} (  1  +  \cot( \theta_m' /2 )   ) 
&   \theta_m' < \pi/2
\\
 2 &
 \theta_m' \geq \pi/2
\end{cases} 
\end{equation}
where $\theta_m' = \min_{x,k,k'} \theta_{m,k,k'}(x)$.

Consequentially, the closure of the normal fibers $  \mathcal{I}_{\mathcal{N}_k(x), \epsilon} $ and $  \mathcal{I}_{\mathcal{N}_{k'}(x'), \epsilon} $ do not touch for two distinct fattened pages $M_{k, S, \epsilon }$ and $M_{k', S , \epsilon}$.
\end{definition}

\begin{proposition} \label{P:star_shaped}
Each cross-section $\omega_{m, \epsilon}(x)$ (as a $2D$ region in the normal plane of $E_m$) is star-shaped with respect to a $2$-ball of radius $c_1 \epsilon$ and contained in a disk of $2$-ball of radius $c_2 \epsilon$ where $c_1$ and $c_2$ are uniform with respect to $x$ and $\epsilon$.
\end{proposition}

\textbf{Proof:} 
This is clear. 
We note $c_1$ depends on the curvature of the pages and the choice of $\epsilon_0$.




\begin{definition} \label{D:special_Lipschitz_domain_I}
A domain $\Omega \subset \mathbb{R}^n$ is called a special Lipschitz domain if there is an orthogonal transformation $T$ of Cartesian coordinates such that
\begin{equation}
 T\Omega = \{ x = (x', x_n) \in \mathbb{R}^n : x' \in \mathbb{R}^{n-1}, \, x_n > \varphi(x') \} 
 \end{equation}
where $\varphi$ is a uniformly Lipschitz function on $\mathbb{R}^{n-1}$.
We call $\varphi$ the boundary graph function to $\Omega$.
\end{definition}

This following proposition follows from our definition of the fattened binding.
The statements in the proposition establish the requirements needed for some embedding and extension theorems.

\begin{figure}[ht!]
\centering
\includegraphics{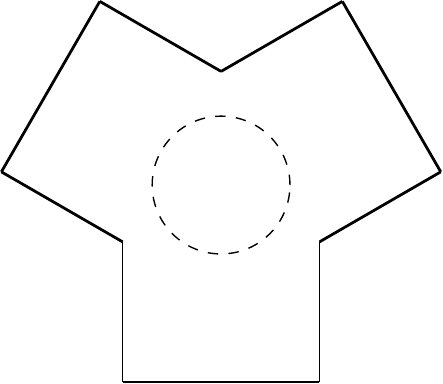}
\caption{A view of $\omega_{m, \epsilon}(x)$ and the ball it is star-shaped with respect to.}
\label{F:starshaped_I}
\end{figure}


\begin{proposition} 
\label{P:Eme_partition_I}
Let $\{E_{m,\epsilon}\}$ be a family of fattened binding neighborhoods as previously described.
For each $\epsilon \in (0, \epsilon_0]$ there exists a partition of unity $\{ \varphi_{i,\epsilon} \}$ ($i \leq N_{U, \epsilon}$, $N_{U, \epsilon}$ depends on $\epsilon$) subordinate to the finite open cover $\{ U_{i,\epsilon} \}$ of $E_{m,\epsilon}$ with the following properties:
\begin{enumerate}
\item $\bigcup_i U_{i,\epsilon}$ is contained in $\bigcup_{x \in E_m} B(x, c_0 \epsilon)$.
\item Each point contained in the covering is in at most $c_N$ sets.
In this sense we say the finite intersection property of these coverings holds in the $\epsilon \to 0$ limit.
\item Each open set $U_{i,\epsilon}$ contains a ball of radius $c_1 \epsilon$ and is contained in a ball of radius $c_2 \epsilon$.
\item If $x \in \partial E_{m,\epsilon}$, then $B(x, c_3 \epsilon ) \subset U_{i,\epsilon}$ for some $i$ and $U_{i,\epsilon} \cap \partial E_{m,\epsilon}$ is a connected subset of some special Lipschitz domain $\Omega_{i, \epsilon}$ whose boundary graph function $\phi_{i, \epsilon}$ has a (Lipschitz) norm bounded above by a constant $c_M$.
\item There is a positive constant $c_\varphi$ such that for each $\epsilon$ the gradient of each $\varphi_{i,\epsilon}$ has a uniform bound $c_\varphi \epsilon^{-1}$:
\begin{equation}
| \nabla \varphi_{i,\epsilon} | \leq c_{\varphi} \epsilon^{-1} .
\end{equation}
\end{enumerate}
\end{proposition}

Most of the properties are obviously compatible and hold for a simple model domain of a homothetically shrinking cylinder.
Point $(4)$ deserves some remark. 
Because of the transversality of the pages, our choice of sleeve length $a_m$, and the fact $M_\epsilon$ is uniformly fattened,
all the ``features'' of the boundary each cross-section $\omega_{m, \epsilon}(x)$ are $O(\epsilon)$ scale (see Figs. \ref{F:Eme_I} and \ref{F:starshaped_I}).
Clearly, this condition would not hold in general if two pages met tangentially.

We will return to this partition of unity later.
The purpose of this partition of unity is to adapt the well-known theorem attributed to Calder\'on and later improved on by Stein \cite{Calderon, Stein} regarding boundedness of extension operators to our shrinking domains.

\begin{theorem} \label{stein}
Let $\Omega$ be an open set in $\mathbb{R}^n$ and let there be positive numbers $r$, $m$, $N$ (an integer) and a sequence $\{ U_i \}_{i \geq 1}$ of open sets satisfying the conditions:
\begin{enumerate}
\item if $x \in \partial \Omega$, then $B(x, r) \subset U_i$ for some $i$,
\item every point $x \in \mathbb{R}^n$ is contained in at most $N$ sets $U_i$,
\item for any $i \geq 1$ there is a special Lipschitz domain $\Omega_i$ with boundary graph function $\varphi_i$ such that $U_i \cap \Omega = U_i \cap \Omega_i$ and
\begin{equation}
|\varphi_i (x') - \varphi_i(y') | \leq m |x' - y'|, \: x', \, y' \in \mathbb{R}^{n-1} 
.
\end{equation}
\end{enumerate}
Then there exists a linear operator $E$ mapping functions defined on $\Omega$ into functions defined on $\mathbb{R}^n$ and having the following properties:
\begin{enumerate}
\item $Eu \, \big|_{\Omega} = u$.
\item $E$ is a continuous operator: $\bigcap_{0 \leq k \leq l} L_p^k ( \Omega) \to \bigcap_{0 \leq k \leq l} L^k_p (\mathbb{R}^n)$ for all $1 \leq p \leq \infty$ and a positive integer $l$ .
\item The norm $||E||_{V^l_p (\Omega) \to V_p^l(\mathbb{R}^n)}$ ($V_p^l(\Omega) := \bigcap_{0 \leq k \leq l} L_p^k ( \Omega)$) is bounded by a constant depending only on $n$, $p$, $l$, $r$, $m$, $N$.
\end{enumerate}
 \end{theorem}

 This theorem has been extended to more general domains \cite{Jones, Rodgers}.
 We approach constructing a family of extension or trace operators by carefully  rescaling each subset of the covering in Proposition \ref{P:Eme_partition_I}.

\subsection{The Fattened Binding Foliation} \label{SS:fbf_I}

\begin{figure}
\centering
\includegraphics{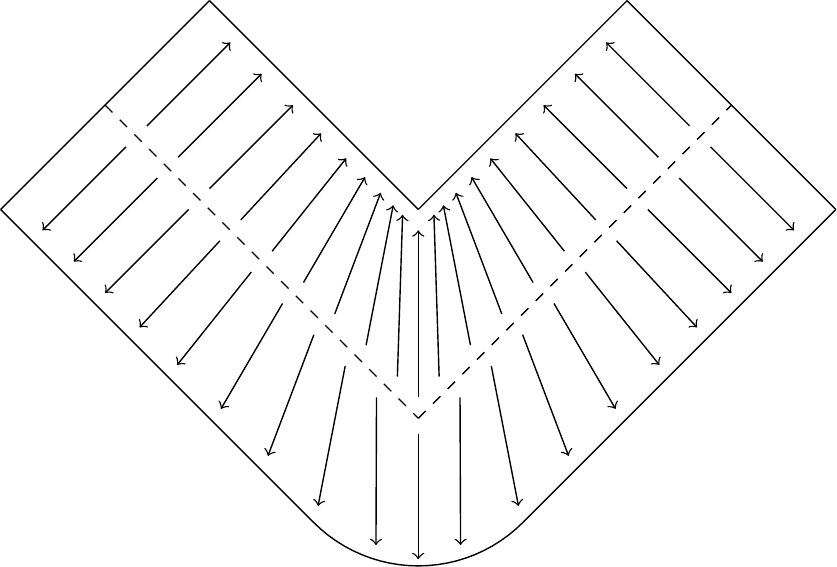}
\caption{Cross sectional view of a pair of vector fields on each of the sleeves yielding a foliation of uniformly fattened binding.}
\label{F:extension_uniform}
\end{figure}

Given our foliations of the fattened pages $M_{k, S, \epsilon}$ and $M_{k', S, \epsilon}$ (in terms of the normal lines $\mathcal{I}_{\mathcal{N}_k(x), \epsilon}$), we wish to extend these foliations into $E_{m, \epsilon}$.
We accomplish this by introducing regions of the fattened binding called sectors.
Breaking up the fattened binding into sectors, we can describe a vector field whose image ``connects'' the foliation of one fattened page to another foliation
 (see Fig. \ref{F:extension_uniform}).

\begin{figure}[ht!]
\centering
\includegraphics{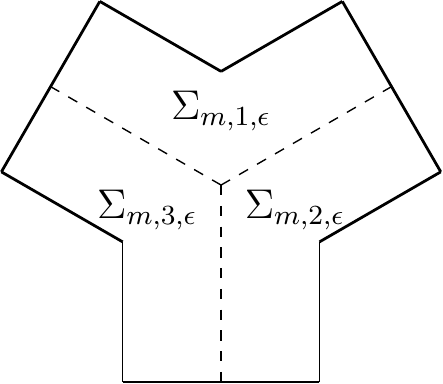}
\caption{Sectors.}\label{F:sect_I}
\end{figure}

\begin{definition} \label{D:sector_I}
Let $E_m$ be a binding and $\{ M_k \}$ ($k \leq n_m$) is the collection of at least two pages that meet at $E_m$ all of which are orientable.
We call the connected components of 
$E_{m, \epsilon} \backslash (E_m \bigcup (\bigcup_{k} S_{k, m, \epsilon}) )$ \textbf{sectors}, and we denote them as $\{ \Sigma_{m, i, \epsilon} \}$ for $i \leq n_m$.
A sector's boundary contains two sleeves of which we say that pair is associated with that sector (see Fig. \ref{F:sect_I}).
\end{definition}

If $E_m$ is a binding connected to non-orientable pages, then taking a partition into local neighborhoods is sufficient for our discussion.
 The case of only one page meeting at a binding is handled separately.

\begin{definition} \label{D:v_I}
Let $E_m$ be a binding and $\{ M_k \}$ ($k \leq n_m$) is the collection pages that meet at $E_m$ all of which are orientable and there are at least two such pages.
We say that the image of family of vector fields $\{t v_{m, i,\epsilon}\}$ ($ t \in (0, 1)$) 
\begin{equation}
v_{m, i, \epsilon} ( x) : E_m \cup S_{k, m, \epsilon} \cup S_{k', m, \epsilon} \mapsto \R^3
\qquad S_{k,m, \epsilon}, \, S_{k', m, \epsilon} \subset \partial \Sigma_{m, i, \epsilon}
\end{equation}
is a \textbf{foliation of the sector matching the foliation of fattened pages} (see Fig. \ref{F:extension_uniform}) if:
\begin{enumerate}
\item $v_{m, i, \epsilon}$ is Lipschitz. 
\item $x \mapsto x + v_{m, i, \epsilon}(x) $ is a homeomorphism between the domain of $v_{m, i, \epsilon}$ and the outward boundary of the sector:
\linebreak
$\partial \Sigma_{m, i, \epsilon} \cap \partial \left( E_{m, \epsilon} \backslash  \bigcup_{k} \partial M_{k, S, \epsilon} \right)$
\item The limit of $v_{m, i, \epsilon}(x)$ as $x \to x' \in \partial S_{k, m, \epsilon} \cap M_k$
is $ \pm \epsilon \mathcal{N}_k(x')$.
\end{enumerate}

If $E_m$ is attached to only one page $M_k$, we say a family of vector fields $\{ v_{m, i, \epsilon } \}$ ($ i = 1, 2$), we impose additionally that the limits of $v_{m, 1, \epsilon}(x)$ and $v_{m, 2, \epsilon}(x)$ match at $E_m$.

\end{definition}

We expand on $(2)$ and describe the construction of functions $\{ v_{m, i, \epsilon} \}$ for all small, positive $\epsilon$ that have uniformly bounded gradients (where they exists).

\begin{proposition} \label{P:v_I}
There is a family of vector-valued functions $\{ v_{m, i, \epsilon} \}$ ($\epsilon \in (0, \epsilon_0 ] $) that extends the foliation of the fattened pages that has length of $O(\epsilon)$ and uniformly bounded gradient (where it exists).
I.e. there exists a $c_1$ and $c_2$ such that
\begin{equation} \label{E:v}
\max_{x \in D(v_{m, i, \epsilon} )} |  v_{m, i, \epsilon}(x) |  \leq c_1 \epsilon,
\end{equation}
and
\begin{equation} \label{E:Dv}
\max_{x \in D(v_{m, i, \epsilon} )} |  \nabla v_{m, i, \epsilon}(x)  | \leq c_2 .
\end{equation}
\end{proposition}
\textbf{Proof:}
We can construct these vector-valued functions following $(2)$ of Definition \ref{D:v_I}.
If $\phi_{m, i, \epsilon}$ is such a homeomorphism per $(2)$, then $v_{m, i, \epsilon}:= \phi_{m, i, \epsilon}(x) - x$ provided the segment from $x$ to $\phi_{m, i, \epsilon}(x)$ is contained in $E_{m, \epsilon}$ (see Fig. \ref{F:extension_uniform}). 
Let us construct such a mapping.

\begin{figure}
\centering
\includegraphics{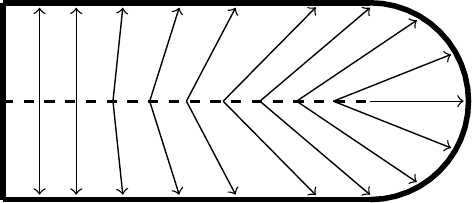}
\caption{Cross sectional view of a pair of vector-valued functions on the sleeves that yield a foliation of fattened binding.}
\label{F:single_page}
\end{figure}

Suppose $\gamma_{1, i, m, x}(s)$ and $\gamma_{2, i, m, x}(t)$ are unit-speed parameterizations of $D(v_{m, i, \epsilon}(x)) := \omega_{m, \epsilon}(x) \cap  (S_{k, m, \epsilon} \cup  S_{k', m, \epsilon}) $ and $R(v_{m, i, \epsilon}(x)) := \omega_{m, \epsilon}(x) \cap \partial \Sigma_{m, i, \epsilon} \cap \partial \left( E_{m, \epsilon} \backslash  \bigcup_{k} \partial M_{k, S, \epsilon} \right)$ respectively.
We denote the length of $R(v_{m, i, \epsilon}(x))$ (as a curve in $\R^3$) as $l \epsilon$.
Then for sufficiently small $\epsilon$, the following mapping has the required segment property ($(x, \phi_{m,i,\epsilon}(x))$ is contained in the sector) for the case of at least two pages meeting at the binding:
\begin{equation}
\phi_{m, i, \epsilon}(
x) =
\gamma_{2, i, m, x} \left( \dfrac{ l  }{ 2 a_{m} }  \gamma_{1, i, m, x}^{-1}(s)  \right)
.
\end{equation}
The case of only one page meeting at the binding requires small modification which one can infer from Fig. \ref{F:single_page}. 
$\Box$

\begin{corollary} \label{C:v_I}
Each sector $\Sigma_{m, i, \epsilon}$ can be parameterized using $v_{m, i, \epsilon}$.
Namely,
a point $x \in \Sigma_{m, i, \epsilon}$ can be written as $x = y + z v_{m, i , \epsilon} (y)$ $(y \in  E_m \bigcup \left(\bigcup_k S_{k, m, \epsilon} \right)$, $z \in (0, 1)$).
\end{corollary}

\subsection{Approximating the Geometry of Fattened Strata}\label{SS:equiv_I}

Here we approximate each fattened page by the product of the corresponding page with an interval.
Although this is not crucial for the proof, we assume that the page $M_k$ is simply connected, otherwise one can partition it further.
Because $M_k$ is partitioned into simply connected patches, the normal $\mathcal{N}_{k}(x)$ is well-defined locally.
A similar analysis is applied to $E_{m}$ and its fattened partner $E_{m, \epsilon}$.

\begin{definition}
\label{D:X_ks_I}
Suppose $U$ is an open region of $\mathbb{R}^2$ with coordinates $y = (y_1,y_2)$.
We define $X_{k, S}$ to be a smooth parameterization of $M_{k, S}$  on $U$:
\begin{equation}
X_{k, S}: (y_1, y_2)\in U\subset \R^2 \mapsto M_{k,S} \subset \mathbb{R}^3 .
\end{equation}
In this subsection, we denote the coefficient functions of the first fundamental form of the immersed surface $M_{k,S}$, (see \cite{Spivak}) as $E$, $F$, and $G$ which are functions on $U$:
\begin{equation}
\label{E:2.40}
\begin{split}
E = D_{y_1} X_{k, S} \cdot D_{y_1} X_{k, S} ,
\\
F = D_{y_1} X_{k, S} \cdot D_{y_2} X_{k, S} ,
\\
G = D_{y_2} X_{k, S} \cdot D_{y_2} X_{k, S} .
\end{split}
\end{equation}
where the symbol ``$\cdot$'' denotes the inner product on $\mathbb{R}^3$.
\end{definition}

\begin{proposition}
\label{P:X_ks_I}
Let $\textbf{e}_{y_1}$ and $\textbf{e}_{y_2}$ denote standard basis vectors of the tangent space $T_{y}U = \R^2$. 
The parameterization $X_{k, S}$ induces a metric $g_{M_{k, S}}$ on $U$.
I.e. $g_{M_{k, S} }$ is the following positive definite bilinear form on $T_y U$: 
\begin{equation} \label{E:g_mks}
g_{M_{k, S} } ( \textbf{a}, \textbf{b} ) := 
\begin{bmatrix} a_1 & a_2 \end{bmatrix} 
\begin{bmatrix} E  &
F
\\
F &
G
\end{bmatrix}
\begin{bmatrix} b_1 \\ b_2 \end{bmatrix}
\end{equation}
where $a_i$, $b_j$ are the respective coefficients of the vectors $\textbf{a}$ and $\textbf{b}$ in the $(\textbf{e}_{y_1}, \textbf{e}_{y_2})$ basis of $T_y U$.
We also use $g_{M_{k, S}}$ to denote the matrix in (\ref{E:g_mks}).
\end{proposition}

\textbf{Proof:}
This is standard (see \cite{Spivak}).

\begin{definition}
\label{D:X_kse}
For sufficiently small $\epsilon$, $M_{k, S, \epsilon}$ admits a parameterization $X_{k, S, \epsilon}$ on $U \times (-\epsilon, \epsilon)$ ($y \in U$, $z \in (-\epsilon, \epsilon)$) where
\begin{equation}
\label{E:2.42}
X_{k, S, \epsilon}(y, z) := X_{k, S}(y) + z \mathcal{N}_{k} ( X_{k,S}(y) ) .
\end{equation}
We denote the coefficient functions of the second fundamental form of an orientable immersed surface, in this case $M_{k, S}$, (see \cite{Spivak}) as $e$, $f$, and $g$:
\begin{equation}
\label{E:2.43}
\begin{split}
 e = - D_{y_1} X_{k, S} \cdot D_{y_1} \mathcal{N}_k ( X_{k, S}( y) ) ,
\\
 f = - D_{y_1} X_{k, S} \cdot D_{y_2} \mathcal{N}_k ( X_{k, S}( y) ) ,
\\
 g = - D_{y_2} X_{k, S} \cdot D_{y_2} \mathcal{N}_k ( X_{k, S}( y) ) .
\end{split}
\end{equation}
\end{definition}

\begin{proposition}
\label{P:X_kse}
The parameterization $X_{k, S, \epsilon}$ induces a metric $g_{M_{k, S, \epsilon}}$ on $U \times (-\epsilon, \epsilon)$:
\begin{equation}
g_{M_{k, S, \epsilon} } :=  
\begin{bmatrix}
E - z e & F - z f & 0  \\
F - z f & G - z g & 0 \\
0 & 0 & 1 
\end{bmatrix}
.
\end{equation}
\end{proposition}

\textbf{Proof:} 
This follows from an explicit calculation of $D_i X_{k, S, \epsilon} \cdot D_j X_{k, S, \epsilon}$ from (\ref{E:2.42}) and simplifying using (\ref{E:2.40}) and (\ref{E:2.43}). $\Box$

\begin{definition} 
\label{D:tildeMke_I}
We define $\tilde{M}_{k,S,\epsilon}$ to be the product space of $M_{k, S}$  and $(-\epsilon, \epsilon)$.
\begin{equation}
\tilde{M}_{k,S,\epsilon} :=  M_{k, S} \times (-\epsilon, \epsilon) .
\end{equation}
\end{definition}

\begin{definition}
\label{D:tildeMke_I}
The product space $\tilde{M}_{k, S, \epsilon}$ admits a parameterization $\tilde{X}_{k, S, \epsilon}$ on an open region $U \times (- \epsilon, \epsilon)$ in $\R^3$ of the form
\begin{equation} \label{E:xkse}
\tilde{X}_{k, S, \epsilon} = ( X_{k, S}, z ).
\end{equation}
\end{definition}

\begin{proposition}
\label{P:tildeMke_I}
The parameterization $\tilde{X}_{k, S, \epsilon}$ (\ref{E:xkse}) induces a metric on $U \times (- \epsilon, \epsilon)$:
\begin{equation}
g_{\tilde{M}_{k, S, \epsilon} } := 
\begin{bmatrix}
E & F & 0 \\
F & G & 0 \\
0 & 0 & 1 
\end{bmatrix}
.
\end{equation}
\end{proposition}

\begin{definition} 
\label{P:phi_Mk_I}
For sufficiently small $\epsilon$, we define a diffeomorphism $\phi_{M_{k, S, \epsilon} }$ from $M_{k,S,\epsilon}$ to $\tilde{M}_{k,S,\epsilon}$ to be given by
\begin{equation}
\phi_{M_{k, S, \epsilon }} (x) =  \tilde{X}_{M_{k, S, \epsilon}} ( X_{M_{k, S, \epsilon}}^{-1} (x)) .
\end{equation}
\end{definition}

\begin{proposition} 
\label{P:tildeMke_I}
The linear operator $\Phi_{M_{k, S, \epsilon } }$ from $H^1(M_{k,S,\epsilon})$ to 
\linebreak
$H^1(\tilde{M}_{k, S, \epsilon})$ induced by the diffeomorphism $\phi_{M_{k, S, \epsilon}}$ (i.e. $\Phi_{M_{k, S, \epsilon} } u =$
\linebreak
$ u (\phi_{M_{k, S, \epsilon } })$) preserves $H^1$-norm of a function up to an $O(\epsilon^{ 1/2})$ error.
\begin{equation} \label{E:a}
\big| \, ||  u ||^2_{H^1(M_{k,S,\epsilon})} - || \Phi_{M_{k, S, \epsilon } } u ||^2_{H^1(\tilde{M}_{k,S,\epsilon})} \, \big| 
\leq c \epsilon  ||  u ||^2_{H^1(M_{k, S,\epsilon})} .
\end{equation}
This inequality (\ref{E:a}) also holds true for other Sobolev spaces $H^n$ and in particular $L_2$.
\end{proposition}

\textbf{Proof:}
First, we show that the metrics $g_{\tilde{M}_{k, S, \epsilon} }$ and $g_{M_{k, S, \epsilon}}$ are close. 

\begin{lemma} 
\label{L:diff_metrics_I}
On the domain $U \times (-\epsilon, \epsilon)$, the two metrics $g_{M_{k, S, \epsilon} }$ and $g_{\tilde{M}_{k, S, \epsilon} }$ are close:
\begin{equation}
 g_{ M_{k, S, \epsilon} } - g_{\tilde{M}_{k, S, \epsilon}}  = B g_{M_{k, S, \epsilon}}
\end{equation}
where matrix $B$ is $O(\epsilon)$ in the Frobenius norm.
\end{lemma}

\textbf{Proof:}
This can be explicitly calculated:
\begin{equation}
\begin{split}
 g_{ M_{k, S, \epsilon} } - g_{\tilde{M}_{k, S, \epsilon}}
=
 \dfrac{ z} {EG - F^2}
 \begin{bmatrix} e   & f     & 0 \\  f & g & 0 \\ 0 & 0 & 0 \end{bmatrix}
  \begin{bmatrix}
 G & -F & 0 \\ -F & E & 0 \\ 0 & 0 & 0
  \end{bmatrix}
g_{\tilde{M}_{k, S, \epsilon}} .
\end{split}
\end{equation}
Because $ | z | \leq \epsilon$, it is clear the right hand side is small.
$\Box$

Having demonstrated the metrics are close, we then calculate the perturbation of two matrix valued functions about $g_{\tilde{M}_{k, S, \epsilon}}$:

\begin{corollary} 
\label{C:diff_metric_I}
The square root of the determinant and inverses of the two metrics $g_{M_{k, S, \epsilon}}$ and $g_{\tilde{M}_{k, S, \epsilon}}$ are also close:
\begin{equation}
\sqrt{\mathsf{det} g_{ {M}_{k, S, \epsilon} } } = \sqrt{ \mathsf{det}g_{  {\tilde{M}}_{k, S, \epsilon} } } \left( 1 + \dfrac{1}{2} \mathsf{det} (   B ) + O(\epsilon^2) \right) ,
\end{equation}
\begin{equation}
{g}_{M_{k, S, \epsilon}}^{-1} = g_{\tilde{M}_{k, S, \epsilon}}^{-1} ( 1 -   B  + O(\epsilon^2)) .
\end{equation}
\end{corollary}

Fixing the target space $U \times (-\epsilon, \epsilon)$ of our coordinate charts on $M_{k,S, \epsilon}$ and $\tilde{M}_{k, S, \epsilon}$, we can now compare functions on $M_{k, S, \epsilon}$ and $\tilde{M}_{k, S, \epsilon}$.

\begin{corollary} 
\label{C:phi_Mk_l2_I}
Let $u \in L_2( U \times (-\epsilon, \epsilon ) )$, then
\begin{align}
\begin{aligned}
\bigg| \int_{ U \times (-\epsilon, \epsilon ) } &|u (y, z) |^2 \sqrt{\mathsf{det}{g}_{M_{k, S, \epsilon}}}\,dydz
\\
&- \int_{U \times (-\epsilon , \epsilon )} | u(y, z) |^2 \sqrt{ \mathsf{det}{g}_{\tilde{M}_{k, S, \epsilon}} }\, dydz \bigg|
\\
&\leq c \epsilon  \int_{U \times ( -\epsilon , \epsilon ) }  | u(y, z) |^2 \sqrt{\mathsf{det}{g}_{M_{k, S, \epsilon}}}\, dydz .
\end{aligned}
\end{align}
\end{corollary}

\begin{corollary} 
\label{C:h1product_I}
Let $u$ be a function in $H^1(U \times (-\epsilon , \epsilon ))$, then
\begin{align}
\begin{aligned}
\bigg| & \int_{U \times (-\epsilon , \epsilon ) } (\nabla u)^* g_{M_{k, S, \epsilon}}^{-1} \nabla u \sqrt{ \mathsf{det} g_{M_{k, S, \epsilon}}} \, dy dz
\\
& - \int_{U \times (-\epsilon , \epsilon ) } (\nabla u)^* g_{\tilde{M}_{k, S, \epsilon}}^{-1} \nabla u \sqrt{\mathsf{det}{g}_{\tilde{M}_{k, S, \epsilon}} } \, dydz \bigg|
\\
& \leq 
c \epsilon  \int_{U \times (-\epsilon , \epsilon ) } (\nabla u)^* g_{M_{k, S, \epsilon}}^{-1} \nabla u \sqrt{ \mathsf{det} g_{M_{k, S, \epsilon}}} \, dy dz
\end{aligned}
\end{align}
where  $\nabla u = ( D_{y_1} u, D_{y_2} u, D_z u)$.
\end{corollary}
These last two statements prove Proposition \ref{P:tildeMke_I}. $\Box$

The cross-sections $\omega_{m, \epsilon}(x)$ vary with $x \in E_m$ due to the curvature of the pages.
Consequentially, more work is needed in defining the parameterization of $E_{m, \epsilon}$.

\begin{definition} 
\label{L:Eme_parameterization_I}
Let $\gamma_m (y)$ be a smooth parameterization of $E_m$ on $U = (0 , l_{E_m} )$.
Let $\textbf{e}_{z_1}$ and $\textbf{e}_{z_2}$ denote standard basis vectors in the normal planes $N_x$ of $E_m$  ($(z_1, z_2) = z \in \omega_{m, \epsilon}(x)$).

We define a fibration $\tilde{U}$ over $U$ as follows:
\begin{equation}
\tilde{U} := \coprod_{ y \in U} \omega_{m, \epsilon}( \gamma_m( y ) )
.
\end{equation}
Let $\Omega_{m, \epsilon} (y , z) :=  z_1 \textbf{e}_{z_1}(\gamma_m(y))  + z_2 \textbf{e}_{z_2}(\gamma_m(y))$.
We can parameterize the $E_{m, \epsilon}$ with $\gamma_m$ and $\Omega_{m, \epsilon}$:
\begin{equation}
Y_{ m, \epsilon}(y, z) = \gamma_m (y) + \Omega_{m, \epsilon}(y, z) .
\end{equation}
\end{definition}


\begin{proposition} 
Parameterization $Y_{m, \epsilon}$ in Definition \ref{L:Eme_parameterization_I} has the following properties:
\begin{enumerate}
\item The image of $\gamma_m (y) + \Omega_{m, \epsilon}(y, \cdot )$ is $\omega_{m, \epsilon}(\gamma_m(y) )$.
\item $D_{z_1} \Omega_{m, \epsilon}(y, \cdot) $ and $D_{z_2} \Omega_{m, \epsilon}(y, \cdot)$ lie in the normal plane of $E_m$ at $\gamma_m(y)$.
\item  $c_1  \leq |   D_{z_i} \Omega_{m, \epsilon}  | \leq c_2$.
\item $ |  D_{y} \Omega_{m, \epsilon} | \leq c_3 \epsilon$.
\end{enumerate}
The parameterization $Y_{ m, \epsilon}$ induces a metric $g_{E_{ m, \epsilon}}$ on $\tilde{U}$:
\begin{equation}
g_{E_{ m , \epsilon } } = \begin{bmatrix} 
1 + D_y \gamma_m \cdot D_y \Omega_{m, \epsilon} & 0 & 0 \\
0 & 1 &
0 
\\
0 & 
0 &
1
\end{bmatrix} 
.
\end{equation}
\end{proposition}

\textbf{Proof}:
Since $\textbf{e}_{z_1}$ and $\textbf{e}_{z_2}$ are orthogonal to $D_y \gamma_m$, $D_y \gamma_m \cdot D_z \Omega_{m, \epsilon} = 0$. 
Because of vectors $\textbf{e}_{z_1}$ and $\textbf{e}_{z_2}$ are orthogonal, we have:
$D_{z_i} \Omega_{m, \epsilon} \cdot D_{z_i} \Omega_{m, \epsilon} = 1$,
 and
$D_{z_i} \Omega_{m, \epsilon} \cdot D_{z_j} \Omega_{m, \epsilon} = 0$ for $i \neq j$. 
$\Box$

\begin{definition} \label{D:tildeEme_I}
We denote by $\tilde{E}_{ m,\epsilon}$ the fibration of $E_{m}$ with fibers $\omega_{m, \epsilon}( x)$:
\begin{equation}
\tilde{E}_{ m, \epsilon} := \coprod_{x \in E_{m} } \omega_{m, \epsilon}(x) .
\end{equation}
\end{definition}

\begin{proposition}
The fibration $\tilde{E}_{m, \epsilon}$ admits a parameterization $\tilde{Y}_{m, \epsilon}$ on $\tilde{U}$:
\begin{equation}
\tilde{Y}_{m, \epsilon} = ( y,  z_1, z_2 ) 
\end{equation}
with an induced metric
\begin{equation}
g_{\tilde{E}_{m, \epsilon} } = \mathsf{Id}_{\R^3}
\end{equation}
where $\mathsf{Id}_{\R^3}$ is the identity matrix.
\end{proposition}

\begin{proposition} \label{P:tildeEme_I}
For sufficiently small $\epsilon$,
there exists a diffeomorphism $\phi_{E_{m, \epsilon } }$ from $E_{m,\epsilon}$ to $\tilde{E}_{m,\epsilon}$ such that the induced linear operator $\Phi_{E_{m, \epsilon } }$ on $H^1(E_{m,\epsilon})$ (i.e. $\Phi_{E_{m, \epsilon } } u = u (\phi_{E_{m, \epsilon } })$) preserves $H^1$-norm up to an $O(\epsilon^{1/2})$ error:
\begin{equation}
\big| \, || u ||^2_{H^1(E_{m,\epsilon})} - || \Phi_{E_{m, \epsilon } } u ||^2_{H^1(\tilde{E}_{m,\epsilon})} \, \big|
\leq c \epsilon  || u ||^2_{H^1(E_{m,\epsilon} )} .
\end{equation}
This inequality also holds true for other Sobolev spaces $H^n$ and in particular $L_2$.
\end{proposition}

\subsection{Bounds on the Sleeves} \label{SS:sleeve_I}

This subsection introduces two needed inequalities.
The proof of the first inequality uses the calculations on induced metrics to show that stretching $M_{k, S}$ back to $M_k$ induces only a small change of a function's norm.
The second inequality involves bounding the $L_2$-norm of a function on a sleeve its $H^1$-norm on the page.

\begin{Pro} \label{P:MkstoMk_I}
There exists a diffeomorphism $\psi_{M_k}$ from $M_k$ to  $M_{k,S}$ such that
\begin{itemize}
\item each column vector of the Jacobian of $\psi_{M_k}$ has length $1 + O(\epsilon)$,
\item for any unit speed differentiable curve $\gamma$ on $\bar{M}_{k,S}$ that is normal to $\partial M_{k, S}$, 
its image
 $\psi_{M_k}(\gamma)$ has unit speed and is normal to the boundary $\partial M_{k}$,
\item the induced operator $\Psi_{M_k}$ (i.e. $\Psi_{M_k} u = u (\psi_{M_k})$) preserves $H^1$-norm up to an $O(\epsilon^{1/2})$ error:
\begin{equation}
\big| \, ||  u ||^2_{H^1(M_k)} - || \Psi_{M_k} u ||^2_{H^1(M_{k,S})} \, \big|  \leq  c \epsilon   || u ||^2_{H^1(M_{k })} .
\end{equation}
This inequality also holds true for other Sobolev spaces $H^n$ and in particular $L_2$.
\end{itemize}
\end{Pro}

\textbf{Proof:}
A sufficiently small neighborhood $V$ of $\partial M_k$ admits a normal coordinate system,
i.e. there is a parameterization $X_k$ on $U \subset \R^2$ of $V$:
\begin{equation} \label{E:X_k}
\begin{split}
X_{k}: (y_1, y_2) \in U =  (0, l_{E_m}) \times  (0, a)  \mapsto M_k ,
\\
\mathrm{such}\,\mathrm{that} \quad  \mathsf{dist}_{M_k}(E_m, X_k(y_1, y_2) ) = y_2 .
\end{split}
\end{equation}
For sufficiently small $\epsilon$, the set $\partial M_{k, S}$ is contained in $V$.
By Definition \ref{D:Ske}, $\partial M_{k,S} \cap E_m$ is the image of $X_k( \cdot, a_m \epsilon)$.
We define a smooth ``shortening'' function
\begin{equation} \label{E:phi_e}
\begin{split}
\varphi_\epsilon:  (0 , a) \mapsto  ( a_m \epsilon, a) 
\quad
\mathrm{such}\,\mathrm{that} \quad
D\varphi_\epsilon \geq 0 ,
\\
 D \varphi_\epsilon ( 0 )  = D \varphi_{\epsilon} (a) = 1 , 
\qquad
| D \varphi_\epsilon - 1 | \leq  c \epsilon
\end{split}
\end{equation}
for some $c>0$.
We can now construct $\psi_{M_k}$:
\begin{equation}
\begin{split}
\psi_{M_k} ( x):= X_k ( (y_1, \varphi_\epsilon (y_1) ) ),
 \quad 
 \mathrm{where} 
 \quad
  ( y_1,  y_2 ) = X_k^{-1} (x) .
\end{split}
\end{equation}
The remainder of the proof follows from the calculating the induced metric from $\psi_{M_k}$ (see Corollaries \ref{C:phi_Mk_l2_I} and \ref{C:h1product_I}). $\Box$

\begin{proposition} \label{P:S_small_I}
Let $M_k$ be a smooth page with boundary $\bigcup_m E_m$.
The $L_2$-norm of a function on $S_{k, m,\epsilon}$ is $O(\epsilon^{1/2})$-bounded by the function's $H^1$-norm on $M_k$:
\begin{equation}
  \int_{S_{k ,m, \epsilon}}  |u |^2  \,    \, dM_k    \leq  c \epsilon^{ }   \int_{M_k }  |u |^2 + |\nabla_{M_k} u|^2 \,   dM_k .
 \end{equation}
\end{proposition}

\textbf{Proof:}
Using the triangle inequality, we conclude
\begin{align}
\begin{aligned}
\int_{S_{k,m,\epsilon}} &|u|^2 \, dM_k
\label{E:psi_mk}
 \leq
  O( \epsilon  )|| u ||_{L_2(M_k)}^2
\\
  &+
 \big| \big| u - \Psi_{M_k } u \big| \big|_{L_2(M_{k, S})}
   \big| \big| u + \Psi_{M_k } u \big| \big|_{L_2(M_{k, S})}
   .
   \end{aligned}
   \end{align}

To bound $ || u - \Psi_{M_k } u ||_{L_2(M_{k, S})}$, we use the coordinate system provided in (\ref{E:X_k}).
Let $X_k$ be the coordinate patch on $U =   (0, l_{E_m}) \times  (0, a)$ (\ref{E:X_k}),
and $\varphi_{\epsilon}$ be the smooth shortening function from (\ref{E:phi_e}).
We define a family of curves that go from $y$ to $\psi_{M_k}(y)$:
\begin{align}
\begin{aligned}
 \gamma_{\varphi_{\epsilon} , y}: & t \in [0, 1] \mapsto U \qquad \gamma_{\varphi_{\epsilon}, y} (0) = ( y_1, y_2)
\\
 &\gamma_{\varphi_{\epsilon} , y} (1) =  (y_1, \varphi_{\epsilon}(y_2)) .
\end{aligned}
\end{align}
 
 In particular we can choose $\gamma_{\varphi_{\epsilon}, y}$ to be constant speed.
Outside of $X_k(U ) \subset M_k$, $ u = \Psi_{M_k} u$, 
and so we need to concern ourselves only with the function on $X_k(U)$.
Let $U' = (0, l_{E_m} ) \times (a_m \epsilon, a)$ and
 let $v(y_1, y_2) = u( X_k (y_1, y_2))$.
 Then we have
\begin{align}
\begin{aligned}
 \big| \big| u -   \Psi_{M_k} u   \big| \big|_{L_2(M_{k, S } ) }^2
= \int_{ U' }  &  \bigg|   \int_{0}^{ 1}   D_{y_2} v \big( y_1, y_2 + t (\varphi_{\epsilon }(y_2) - y_2)   \big) 
\\
&| \varphi_{\epsilon}(y_2 ) - y_2 |  dt \bigg|^2 
 \sqrt{\mathsf{det}{g_{M_k}} }  dy
.
\end{aligned}
\end{align}
Let $\xi =  y_2 + t (\varphi_{\epsilon} ( y_2 ) - y_2 )$,
and so $ d\xi = dy_2 ( 1 - t + t  D \varphi_{\epsilon} (y_2) ) $.
Because $D \varphi_{\epsilon} (y_2) = 1 + O(\epsilon )$,
we can then write $d\xi = dy_2( 1 - t O(\epsilon  ))$.
Thus, the Jacobian $J$ from $(y_1, y_2) \mapsto (y_1, \xi)$
is of the form $1 + O(\epsilon)$.
Applying $| \varphi_{\epsilon } (y_2) - y_2 | = O(\epsilon  )$, we have
\begin{align} \label{E:u.66}
\begin{aligned}
\big| \big| u  & - \Psi_{M_k} u \big| \big|_{L_2(M_{k, S } ) }^2
\\
 & \leq
  \int_{ U'  } \,    \int_{0}^{ 1}   \Big| D_{\xi} v (y_1, \xi )   \Big|^2
 \dfrac{ O(\epsilon^2 ) }{1 - t O(\epsilon  )} \, dt \, \sqrt{\mathsf{det}{g_{M_k}} }(y_1, \xi )  \, dy_1 d\xi
\\
&  \leq 
   \big| \big| D_{\xi }  v  \big| \big|^2_{L_2( U' , \mathsf{det}g_{M_k}^{1/2} ) } \int_{0}^{ 1}
\dfrac{O(\epsilon^2  ) }{1 - t O( \epsilon  )} \, dt
\\
&\leq 
O(\epsilon^2) || \nabla_{M_k} u ||^2_{L_2(M_k) } .
\end{aligned}
\end{align}
Applying (\ref{E:u.66}) to (\ref{E:psi_mk}) we get the $O(\epsilon)$ bounds.
$\Box$

\subsection{Local Extensions of Functions on a Stratum to the Fattened Domain} \label{SS:ext_I}

We can extend a function form $M_{k,S}$ into $M_{k, S, \epsilon}$ by first extending along the fibers and then applying the diffeomorphism operator in Proposition \ref{P:tildeMke_I}.
The extension from the binding and the sleeves is handled by extending along the foliation derived in Definition \ref{D:v_I} by means of its associated coordinate system (Corollary \ref{C:v_I}).

\begin{definition} \label{D:tilde_EzMkse_I}
We define $\tilde{\mathcal{E}}_{k,z,\epsilon}$ to be the extension operator  from $M_{k,S}$ to $\tilde{M}_{k,S,\epsilon}$, a bounded linear operator from $L_2(M_{k, S})$ to $L_2( \tilde{M}_{k, S, \epsilon})$,
given by:
\begin{equation} 
\mathcal{E}_{k,z,\epsilon} u(y,z) = u(y)
.
\end{equation}
where $y \in M_{k,S}$  and $z \in ( - \epsilon, \epsilon)$.
\end{definition}

\begin{definition} \label{D:EzMkse_I}
Let $u \in L_2(M_{k,S})$.
We define ${\mathcal{E}}_{k,z,\epsilon}$ to be the extension operator from $M_{k,S}$ to ${M}_{k,S,\epsilon}$
given by 
\begin{equation}
\mathcal{E}_{k,z,\epsilon} := \Phi^{-1}_{M_{k, S, \epsilon}} \tilde{\mathcal{E}}_{k,z, \epsilon} .
\end{equation}
\end{definition}



\begin{proposition} \label{P:E_kze_I}	
For $u \in H^1(M_{k,S} )$, one has:
	\begin{equation}	
		\big| \, || u ||^2_{L_2(M_{k,S} )} - || (2\epsilon)^{-1/2} {\mathcal{E}}_{k, z, \epsilon} u ||^2_{ L_2(
 {M}_{k,S,\epsilon} ) }\, \big|
\leq c \epsilon ||u ||^2_{L_2(M_{k, S})}
	\end{equation}
	and
		\begin{equation}	
		\big| \, || \nabla_{M_k} u ||^2_{L_2(M_{k,S} )} - || \nabla (2\epsilon)^{-1/2}  {\mathcal{E}}_{k, z, \epsilon}  u ||^2_{ L_2(
 {M}_{k,S,\epsilon} ) }\, \big|
\leq c \epsilon || u ||^2_{H^1(M_{k, S})} .
	\end{equation}
\end{proposition}

\textbf{Proof:}
Because
\begin{equation}
\int_{-\epsilon}^{\epsilon} \dfrac{1}{2 \epsilon} | \tilde{\mathcal{E}}_{k, z, \epsilon} u(y, z)|^2 \, d z = | u(y) |^2
,
\end{equation}
it follows that
\begin{align}
\label{E:68}
\begin{aligned}
		\big| \, || u & ||^2_{L_2(M_{k,S} )}- || (2 \epsilon)^{-1/2} \tilde{\mathcal{E}}_{k, z, \epsilon} u ||^2_{ L_2(
\tilde{M}_{k,S,\epsilon} ) }\, \big|  
\\
& 
=
\bigg| \,  \int_{M_{k,S}}  | u |^2 \,   dM_k - \int_{ {M}_{k , S } } \int_{  -\epsilon  }^{\epsilon}  \dfrac{1}{2\epsilon } | \tilde{\mathcal{E}}_{k, z, \epsilon} u(y, z, \epsilon ) |^2 \, dz \, d{M}_{k} \, \bigg|  
\\
& 
= 0 
.
\end{aligned}
\end{align}
Turning to the norm of the gradient, we have	
	\begin{align} 
	\label{E:bbb}
	\begin{aligned}
		&\big| \, || \nabla_{M_k} u ||^2_{L_2(M_{k,S} )} - || \nabla (2\epsilon)^{-1/2} \tilde{\mathcal{E}}_{k,z, \epsilon} u ||^2_{ L_2(
\tilde{M}_{k,S,\epsilon} ) }\, \big|
\\
  &=
\bigg| \,  \int_{M_{k,S}}  | \nabla_{M_k} u |^2 \,    dM_k -      \int_{ {M}_{k , S } } \int_{   -\epsilon  }^{\epsilon }   \dfrac{1}{2\epsilon } | \nabla \tilde{\mathcal{E}}_{k, z, \epsilon } u |^2 \, dz \, d\tilde{M}_{k } \, \bigg|   
.
\end{aligned}
\end{align}
Clearly $D_z \tilde{\mathcal{E}}_{k, z, \epsilon} u = 0$, so we can rewrite (\ref{E:bbb}) to get:
\begin{align}
\label{E:71}
\begin{aligned}
\bigg| \,  \int_{M_{k,S}}  | \nabla_{M_k} &u |^2 \,    dM_k
\\
&-
\dfrac{1}{2\epsilon}  \int_{ {M}_{ k, S } }     \int_{-\epsilon  }^{ \epsilon }  | \nabla_{M_{ k }} \tilde{\mathcal{E}}_{k, z, \epsilon} u |^2 + | D_z \tilde{\mathcal{E}}_{k, z, \epsilon } u |^2 \, dz \, d\tilde{M}_{k } \, \bigg|
\\
&= 0 
. 
\end{aligned}
\end{align}
Then we apply Proposition \ref{P:tildeMke_I} to the results in (\ref{E:68}) and (\ref{E:71}).
$\Box$

\begin{definition} \label{D:EzSkme_I}
We define $\mathcal{E}_{m, S, z, \epsilon}$ to be the extension operator on $L_2(E_m \bigcup \left( \bigcup_k S_{k, m, \epsilon} \right) )$ to $L_2({E}_{m, \epsilon} ) $
given by sector as
\begin{align}
\begin{aligned}
\mathcal{E}_{m, S, z, \epsilon} &u (y, z) 
= 
u(y) 
\end{aligned}
\end{align}
where $(y, z)$ is the coordinate system described in Corollary \ref{C:v_I}.
\end{definition}

 \begin{proposition} \label{P:Ezbdd_I}
 The extension operators $(2\epsilon)^{-1/2} \mathcal{E}_{m, S, z, \epsilon}$ from 
 \linebreak
 $H^1(E_m \bigcup \left( \bigcup_k  S_{k,m,\epsilon} \right) )$ to 
  $H^1( {E}_{m,\epsilon} )$ satisfy the following bound:
 \begin{equation} \label{E:EmSze_I}
 || (2 \epsilon )^{-1/2} \mathcal{E}_{m, S, z,  \epsilon} u ||^2_{H^1( {E}_{m,\epsilon})} 
 \leq c ||u||^2_{H^1(E_m \bigcup ( \bigcup_k  {S}_{k, m, \epsilon} ) )} .
 \end{equation}
 \end{proposition}

 \textbf{Proof:}
 Let $x = y +  z v_{m, i,\epsilon}(y) $ 
 ($ y \in E_m \bigcup$ $\left( \bigcup_k S_{k, m, \epsilon} \right) $, $z \in (0, 1)$).
 We break the sector into two sets, the two halves each ``above'' a sleeve:
  $\{ y + z v_{m, i, \epsilon}(y): y \in S_{k, m, \epsilon}, \;z \in (0,1) \}$ for some page index $k$. 
 We calculate the induced metric on this region to demonstrate the determinant of the metric is the correct order of $\epsilon$ such that the $L_2$ part of ($\ref{E:EmSze_I}$) holds.
To accomplish that, we use the parameterization of $S_{k,m,\epsilon}$ in (\ref{E:X_k}) renaming the parameterized variable as $t$ ($X_k(t) = y$), and we denote the induced metric on the domain of $t$ for $S_{k, m, \epsilon}$ as $g_{M_k}$.
 The induced metric $g_{\Sigma_{m, i, \epsilon}, k}$ on $(t,z)$ is
 \begin{equation}
 g_{\Sigma_{m, i, \epsilon}, k} 
 =
 \begin{bmatrix}  
 g_{M_{k}} + D_t X_{k} \cdot z D_t v_{m, i, \epsilon}(X_k) &   D_t X_k \cdot v_{m, i, \epsilon}
 \\
 D_t X_k \cdot v_{m, i, \epsilon} 
 &
 | v_{m, i, \epsilon} |^2
 \end{bmatrix}
 .
 \end{equation} 
 It follows,
 \begin{equation}
 \mathsf{det}( g_{\Sigma_{m, i, \epsilon}, k} ) 
 \leq c \epsilon \mathsf{det}( g_{M_k})
 .
 \end{equation}
 Thus we conclude
  \begin{equation}
  \label{E:2.83}
  || (2\epsilon)^{-1/2} \mathcal{E}_{m, S, z,  \epsilon} u ||^2_{L_2( {E}_{m,\epsilon})}
   \leq c ||u||^2_{L_2(E_m \bigcup ( \bigcup_k  {S}_{k, m, \epsilon}) ) } .
 \end{equation}
To calculate the gradient at the point,
we calculate the divided difference between $\mathcal{E}_{m, S, z, \epsilon}  u(x )$ and $\mathcal{E}_{ m, S, z, \epsilon} u (x + \delta   )$.
\begin{align}
\begin{aligned}
|\nabla_{\hat{\delta}} \mathcal{E}_{m, S, z , \epsilon}|
&=
\big| \limsup_{\delta \to 0} \dfrac{ \mathcal{E}_{m, S, z, \epsilon} u ( y + \delta_y, z + \delta_z ) - \mathcal{E}_{ m , S, z, \epsilon} u (y , z     )     }{ | \delta |}  \big|
\\
&= \big|  \limsup_{\delta \to 0} \dfrac{ u(y + \delta_y ) - u(y)   }{ |  \delta_y   | } \dfrac{| \delta_y   | }{ | \delta |}  \big|
\leq c_{\nabla v_{m, i, \epsilon}} | \nabla_{M_k} u | .
\end{aligned}
\end{align}
This lets us conclude
\begin{equation}
| \nabla \mathcal{E}_{m, S, z, \epsilon} u(y, z) | \leq  c_{\nabla v_{m, i,\epsilon} } |  \nabla_{M_k} u(y) | .
\end{equation}
Hence we arrive at a bound on the derivative giving us (\ref{E:EmSze_I}) with (\ref{E:2.83}):
\begin{equation}
 || \nabla (2 \epsilon )^{-1/2} \mathcal{E}_{m, S, z,  \epsilon} u ||^2_{L_2( {E}_{m,\epsilon})} \leq \sum_{k} c || \nabla_{M_k}u||^2_{L_2( {S}_{k, m, \epsilon})} . \qquad   \Box
\end{equation}

\subsection{Extension Operator $K_\epsilon$}\label{SS:exp_contr}

Now we can define the extension operators in the sense of Definition \ref{D:Ke}.

\begin{Pro} \label{P:Ke_1_I}
Let $M$ be an open book structure.
Let $\Lambda  \leq c \epsilon^{-1 + \delta}$ where $\delta > 0$ and $\Lambda \notin \sigma (A) $.
For some $\epsilon_0 > 0$, the family of linear operators $\{K_{\epsilon}  \}_{\epsilon \in (0, \epsilon_0 ] }$ that satisfies the conditions in Definition \ref{D:Ke} is ($u \in \mathcal{P}_\Lambda L_2(M)$)
\begin{equation} \label{E:u.83}
K_\epsilon  u := \begin{cases}
(2 \epsilon)^{ -1 / 2 }   \mathcal{E}_{k,z,\epsilon}   u  & M_{k,S}
\\
(2 \epsilon)^{ -1 / 2 }    \mathcal{E}_{ m, S, z, \epsilon} u & E_m \bigcup \left( \bigcup_k S_{k, m, \epsilon} \right) .
\end{cases}
\end{equation}
\end{Pro}

\textbf{Proof:}
Beginning with $E_m \bigcup \left( \bigcup_k S_{k, m, \epsilon} \right)$,
we apply Proposition \ref{P:Ezbdd_I} to get
\begin{equation} \label{E:2.78}
|| ( 2\epsilon) ^{ -1 / 2 }    \mathcal{E}_{ m, S, z, \epsilon} u||^2_{H^1(E_{m, \epsilon})}
\leq 
c || u ||^2_{H^1(E_m \bigcup ( \bigcup_k S_{k, m, \epsilon} ) ) } .
\end{equation}
Applying the spectral embedding Proposition \ref{P:proj}, the previously expression is bounded by
$c (1 + \Lambda) || u ||^2_{L_2(E_m \bigcup ( \bigcup_k S_{k, m, \epsilon} ) ) }$
which in turn is bounded by the energy on $M$ (Proposition \ref{P:S_small_I}). 
This yields an upper bound of
\begin{equation}
c (1 + \Lambda) \epsilon^{ } || u ||^2_{ \mathcal{G}^1  } = o(1) || u ||^2_{ \mathcal{G}^1  } .
\end{equation}
Therefore the right-hand term of (\ref{E:2.78}) is small.
For a function on $M_{k, S}$, we show that $u$ is not only close to its extension $(2 \epsilon)^{-1/2} \mathcal{E}_{k, z, \epsilon} u$ (\ref{E:u.83}) in $L_2$ but also in $H^1$.
Starting with the following norm difference
\begin{equation} \label{E:K_i}
\big|
\sum_k || (2 \epsilon)^{ -1 / 2 }  \mathcal{E}_{k, z, \epsilon}
  u  ||_{H^1(M_{k, S,  \epsilon} ) }
- ||u ||_{ \mathcal{G}^1 } \big| ,
\end{equation}
we break $||u ||_{\mathcal{G}^1}$ into page terms and sleeve terms and use the triangle inequality.
We get an upper bound of (\ref{E:K_i}) of
\begin{align}
\label{E:2.82}
\begin{aligned}
\sum_k \big|
|| (2 \epsilon)^{ -1 / 2 } & \mathcal{E}_{k, z, \epsilon}    u  ||_{H^1(M_{k, S, \epsilon} ) }
- ||u ||_{H^1(M_{k, S} ) } \big|
\\
& +  || u ||_{H^1(E_m \bigcup ( \bigcup_k S_{k, m, \epsilon} ) ) } .
\end{aligned}
\end{align}
The first term of (\ref{E:2.82}) is bounded by Proposition \ref{P:E_kze_I}.
After a norm bound on the sleeve (Propositions \ref{P:S_small_I} and \ref{P:proj}), we conclude (\ref{E:K_i}) is bounded by 
$(1 + \Lambda)^{1/2} O(\epsilon^{ 1/2} )  || u ||_{\mathcal{G}^1 }$.
 Thus $K_\epsilon$ is a near isometry in both $L_2$ and $H^1$.
 $\Box$

\subsection{Local Averaging Operators} \label{SS:ir_I}

This subsection concerns an averaging operator on the fattened page and an averaging operation on the fattened binding constructed by means of an integral representation.
These averaging operators satisfy some Poincar\'e-type inequalities.
I.e. the norm of the difference between a function and a constant (in the simplest formulation this constant is the average) is bounded by the norm of the function's derivative.

\begin{definition} \label{D:N_I}
Let $\tilde{N}_{k, \epsilon} $ denote the following bounded linear operator on $L_2(\tilde{M}_{k, S,\epsilon})$:
\begin{equation}
\tilde{N}_{k , \epsilon } u(y,z) := \dfrac{1}{ 2 \epsilon } \int_{ -\epsilon}^{ \epsilon } u(y, \zeta)  \, d\zeta
\qquad y \in M_{k, S}, \; z \in ( - \epsilon, \epsilon ) .
\end{equation}
\end{definition}

\begin{definition} \label{R:N}
The averaging operator $N_{k, \epsilon}$ on $M_{k,S,\epsilon}$ is given by composition with the diffeomorphism $\phi_{M_{k, S, \epsilon}}$: 
\begin{equation}
N_{k, \epsilon} := \Phi_{M_{k, S, \epsilon } }^{-1} \tilde{N}_{k, \epsilon} \Phi_{M_{k , S, \epsilon} } .
\end{equation}
We also let $N_{k, \epsilon}$ denote a bounded linear operator from $L_2(M_{k, S, \epsilon})$ to $L_2(M_{k, S})$ by restricting $N_{k, \epsilon} u$ to $M_{k, S}$.
\end{definition}

\begin{proposition}
\label{P:N_I}
The norms of the family of averaging operators $\{  {N}_{k, \epsilon} \}$ on $L_2( {M}_{k,S,\epsilon})$ has a uniform bound $c$.
\end{proposition}

\textbf{Proof:}
Boundedness is clear from the Cauchy-Schwartz Inequality.

\begin{proposition} \label{P:Poincare_interval}
For $u \in H^1(\tilde{M}_{k,S,\epsilon})$,
$\tilde{N}_{k, \epsilon}$ satisfies a Poincar\'e-type inequality:
\begin{equation}
\int_{\tilde{M}_{k,S,\epsilon}} | u - \tilde{N}_{k, \epsilon} u |^2 \, d\tilde{M}_{k,S,\epsilon} \leq c \epsilon^2 \int_{\tilde{M}_{k,S,\epsilon}} | \nabla u |^2 \, d\tilde{M}_{k,S,\epsilon} .
\end{equation} 
\end{proposition}

\textbf{Proof:}
Because the lowest non-constant Neumann eigenfunction for the interval $(-1, 1)$ is $\sin(\pi x / 2)$,
the Poincar\'e inequality for an $\epsilon$-interval yields
\begin{equation} \label{E:2.86}
\begin{split}
 \int_{ -\epsilon }^{\epsilon }  | u - \tilde{N}_{ k, \epsilon } u |^2 \, dz  \leq
  \dfrac{4\epsilon^2}{  \pi^2 }    \int_{- \epsilon  }^{\epsilon }  | D_z  u(y, z) |^2 \, dz .
\end{split}
\end{equation}
We then integrate (\ref{E:2.86}) over $M_{k, S}$. 
Because $\tilde{M}_{k, S, \epsilon}$ is a product of $M_{k, S}$ and $( -\epsilon, \epsilon)$, the result follows from Fubini's theorem. $\Box$

\begin{corollary} \label{C:N_Poincare_I}
For $u \in H^1( M_{ k, S, \epsilon} )$, the averaging operator $N_{ k, \epsilon}$ admits a Poincar\'e-type inequality:
\begin{equation} 
\label{E:2.87}
|| u -  N_{ k, \epsilon} u ||^2_{L_2( M_{ k, S, \epsilon} ) } \leq
c \epsilon^2 || \nabla u ||^2_{L_2( M_{ k, S, \epsilon} )} .
\end{equation}
\end{corollary}

 \textbf{Proof:} 
The inequality (\ref{E:2.87}) is straightforward application of Proposition \ref{P:tildeMke_I}.


\begin{proposition} \label{P:Nu_I}
For $u \in H^1({M}_{k,S,\epsilon})$, one has:
\begin{equation}
\bigg| \,  |  | u  |  |_{L_2({M}_{k,S,\epsilon})}^2 
-  |  | (2 \epsilon)^{1/2} {N}_{k, \epsilon} u  |  |_{L_2(M_{k,S} )}^2 \, \bigg| 
\leq c \epsilon 
  |  | u  |  |_{H^1({M}_{k, S,\epsilon})}^2 .
\end{equation}
\end{proposition}

\textbf{Proof:} Bounding the difference squared in the fibered space, we get:
\begin{align}
\begin{aligned}
 |  | u  |  &|_{L_2(\tilde{M}_{k,S,\epsilon})}^2 
-  |  | (2 \epsilon)^{1/2} {N}_{k, \epsilon} u  |  |_{L_2(M_{k,S} )}^2 
\\
&=
\int_{\tilde{M}_{k,S,\epsilon}} |u|^2  \, d\tilde{M}_{k,S,\epsilon}
-  \int_{M_{k, S}} \Big( \int_{ -\epsilon  }^{\epsilon } | \tilde{N}_{ k, \epsilon } u|^2 \, dz \Big) \, dM_k
\\
&\leq (1+ O(\epsilon)) | | u - \tilde{N}_{k, \epsilon } u | |_{L_2(\tilde{M}_{k,S,\epsilon})} 
| | u +  \tilde{N}_{k, \epsilon } u | |_{L_2(\tilde{M}_{k,S,\epsilon})}
\\
&\leq  2 \epsilon (1 + O(\epsilon)) || u||_{H^1(\tilde{M}_{k, S, \epsilon} ) }^2 .
\qquad \Box
\end{aligned}
\end{align}

Then we apply Proposition \ref{P:tildeMke_I} on the above.

\begin{proposition} \label{P:QNu_I}
The linear operator ${N}_{k, \epsilon}$ is bounded on $H^1({M}_{k,S,\epsilon})$,
\begin{equation} \label{E:u.97}
\int_{M_{k,S}} |\nabla_{M_k} (2 \epsilon)^{1/2} {N}_{k , \epsilon } u|^2 \,   dM_{k}
 \leq
 (1 + O(\epsilon))
  \int_{ {M}_{k,S,\epsilon}} | \nabla_{M_k} u |^2 \, d {M}_{ \epsilon} .
\end{equation}
\end{proposition}

\textbf{Proof:}
We begin with demonstrating (\ref{E:u.97}) for a function on the space $\tilde{M}_{k, S, \epsilon}$:
\begin{equation}
\int_{M_{k,S}} |\nabla_{M_k} (2 \epsilon)^{1/2} \tilde{N}_{k, \epsilon} u|^2    dM_{k}  
= \int_{M_{k,S}} \Big( \int_{ - \epsilon }^{ \epsilon }  | \nabla_{M_k} \tilde{N}_{k, \epsilon } u|^2 dz \Big)  dM_k .
\end{equation}
Because
\begin{align} \label{E:u.99}
\begin{aligned}
\int_{\tilde{M}_{k,S,\epsilon}} \Big| D_{y_i} \dfrac{1}{2 \epsilon} \int_{ -\epsilon }^{\epsilon }   u  d \zeta  \Big|^2  d\tilde{M}_{k,S,\epsilon}
= \int_{\tilde{M}_{k,S,\epsilon}} \Big| \dfrac{1}{2 \epsilon} \int_{ -\epsilon }^{ \epsilon }   D_{y_i} u   d\zeta \Big|^2  d\tilde{M}_{k,S,\epsilon},
\end{aligned}
\end{align}
we then use the embedding of $L_1$ in $L_2$ on a compact interval and the Cauchy-Schwartz Inequality:
\begin{align}
\begin{aligned}
&\int_{\tilde{M}_{k, S,\epsilon}} \Big| \dfrac{1}{2 \epsilon} \int_{ -\epsilon }^{ \epsilon }    \nabla_{M_k} u \, d\zeta \Big|^2  \, d\tilde{M}_{k, S, \epsilon}
\\
&\leq
\int_{\tilde{M}_{k,S,\epsilon}}\Big( \dfrac{1}{2 \epsilon} \int_{ -\epsilon }^{ \epsilon } |   \nabla_{M_k} u |^2 \, d\zeta  \Big)  \, d\tilde{M}_{k,S,\epsilon}
\\
&\leq
\int_{\tilde{M}_{k,S,\epsilon}} |\nabla_{M_k} u|^2 \, d\tilde{M}_{k,S,\epsilon} .
\end{aligned}
\end{align}

Then we get the result after an application of Proposition \ref{P:tildeMke_I}. 
$\Box$

\begin{lemma} \label{L:Riesz}
Let $\Omega$ be a bounded domain in $\mathbb{R}^n$ with diameter $D$.
Suppose $l > 0$ and
\begin{equation}
\label{E:2.95}
 Ru(z) := \int_{\Omega} \dfrac{ u(\zeta) }{ |z - \zeta|^{n - l} } d\zeta .
\end{equation}
Then $R$ is a continuous linear operator on $L_p (\Omega)$, $1 \leq p \leq \infty$, and
\begin{equation} \label{E:R}
 ||R|| \leq n \, |B(0,1)| D^l/l 
 .
\end{equation}
\end{lemma}

\textbf{Proof}:
Let $\chi$ be the characteristic function of $B(0, D)$.
Letting our test function be zero outside of $\Omega$ and $K = |z|^{l-n} \chi(z)$, we observe
$Ru(z) = (K \ast u)|_{\Omega}$. 
Therefore the inequality (\ref{E:R}) follows from the Young inequality.
$\Box$

The kernel in (\ref{E:2.95}) appears in the remainder term in the following integral representation (see \cite{MaPo}):

\begin{theorem} \label{T:integral_rep}
Let $\Omega$ be a bounded domain star-shaped with respect to a ball $B(0,\delta) \subset \Omega$ in $\mathbb{R}^n$ and let $u \in L^l_p(\Omega)$.
Then for almost all $x \in \Omega$
\begin{align} 
\label{E:ir}
\begin{aligned}
 u(z) = &\delta^{-n} \sum_{|\alpha| < l} \left( \dfrac{z}{\delta} \right)^{\alpha} \int_{B(0,\delta)} \phi_{\alpha} \left( \dfrac{\zeta}{\delta} \right) u(\zeta) d\zeta  
 \\
 &+ \sum_{|\alpha| = l} \int_{\Omega} \dfrac{ f_{\alpha} (z, r, \theta) }{r^{n-l}} D^{\alpha} u(\zeta) d\zeta
 ,
\end{aligned}
\end{align}
where $r = |z-\zeta|$, $\theta = (\zeta -z)/r$, $\phi_{\alpha} \in C^\infty_0(B(0,1))$, and $f_\alpha$ are infinitely differentiable functions such that
\begin{equation}
 |f_\alpha | \leq c( D/ \delta)^{n-1} .
\end{equation}
$c$ is a constant independent of $\Omega$ and $D$ is the diameter of $\Omega$.
\end{theorem}

\begin{remark} \label{R:rep}
Let $\varphi \in C_0^{\infty}(B(0,1))$ such that $\int_{B(0,1)} \varphi = 1$.
The function $f_\alpha$ in the integral representation (\ref{E:ir}) has an explicit expression in terms of $\varphi$; in particular (\ref{E:ir}) can be written as:
\begin{align}
\begin{aligned}
& u(z) 
 = \delta^{-n} \sum_{|\alpha| < l} \dfrac{1}{\alpha !} \int_{B(0,\delta)} \varphi  \left( \dfrac{\zeta}{\delta} \right) (z - \zeta )^\alpha D^{\alpha} u(\zeta) d\zeta
\\
&+ 
\sum_{|\alpha| = l} \dfrac{(-1)^l l}{\alpha!} \int_{\Omega}  \Big(  \int_{r}^{\infty} \varphi \left( \dfrac{z + \rho \theta}{\delta} \right) \rho^{n-1} \, d \rho \Big)  \dfrac{  \theta^{\alpha} }{r^{n-l}} D^{\alpha} u(\zeta) d\zeta .
\end{aligned}
\end{align}
\end{remark}

\textbf{Proof:}
These are standard results in the theory of differentiable functions \cite{ Ma, MaPo}.

We note this representation (\ref{E:ir}) in particular holds on almost every slice of a fibration like $\tilde{E}_{m,\epsilon}$.

\begin{definition} \label{D:P}
We define $\tilde{P}_{ m , \epsilon }$ to denote the following bounded linear operator on $L_2(\tilde{E}_{m,\epsilon})$.
\begin{equation}
\tilde{P}_{ m , \epsilon }u(y, z ) := \dfrac{1}{|D( 0, c_r  \epsilon )|} \int_{D( 0, c_r  \epsilon )} \varphi \left(\dfrac{ \zeta }{ c_r  \epsilon  }\right)  u(y, \zeta) \, d\omega_{m, \epsilon}(y) \; ,
\end{equation}
where $\varphi \in C_0^{\infty}(D(0,1))$ such that $ \displaystyle  \int_{D(0,1)} \varphi = 1$.
\end{definition}

\begin{definition} \label{R:N}
The averaging operator $P_{m, \epsilon}$ on $E_{ m,\epsilon}$ is given by composition with the corresponding diffeomorphism: 
\begin{equation}
P_{m, \epsilon} := \Phi_{E_{m, \epsilon } }^{-1} \tilde{P}_{m, \epsilon} \Phi_{E_{m, \epsilon} } .
\end{equation}
We also let $P_{m, \epsilon}$ denote a bounded linear operator from $L_2(E_{m, \epsilon})$ to $L_2(E_m)$ by restricting $P_{m, \epsilon} u$ to $E_{m}$.
\end{definition}

\begin{proposition}
The norms of the family of averaging operators $\{{P}_{m, \epsilon} \}$ on $L_2({E}_{ m,\epsilon})$ have a uniform upper bound $c$.
\end{proposition}

As with the operator ${N}_{k, \epsilon} $, boundedness of ${P}_{m, \epsilon}$ is clear from the Cauchy-Schwartz Inequality and Proposition \ref{P:tildeEme_I}.

\begin{Pro} \label{P:QPu}
The linear operator ${P}_{m, \epsilon}$ is bounded on $H^1({E}_{ m,\epsilon})$:
\begin{equation}
\begin{split}
\int_{ {E}_{m, \epsilon} } |\nabla_{   }   {P}_{m, \epsilon}  u|^2 \,   d{E}_{m, \epsilon }
 \leq 
 (1 + O(\epsilon) )  || \varphi ||^2_{B(0, 1)} \int_{ {E}_{ m,  \epsilon}} | \nabla_{  } u |^2 \, d  {E}_{ m ,\epsilon} .
\end{split}
\end{equation}
\end{Pro}

\textbf{Proof:}
We begin with bounding the norm of the derivative of a function on $\tilde{E}_{m, \epsilon}$.
We have 
\begin{align}
\begin{aligned}
&\int_{\tilde{E}_{ m,\epsilon}} \Big| \nabla_{  } \dfrac{1}{| D( 0, c_r  \epsilon ) | } \int_{D(0, c_r \epsilon) } \varphi \left( \dfrac{\zeta}{ c_r \epsilon  } \right) u \, d \omega_{ m, \epsilon} (y)  \Big|^2 \, d\tilde{E}_{ m , \epsilon}
\\
&= \int_{\tilde{E}_{ m,  \epsilon}} \Big| \dfrac{1}{| D(0, c_r \epsilon ) |}
\int_{D(0, c_r \epsilon) } \varphi \left( \dfrac{\zeta}{ c_r \epsilon   } \right) D_y u  \, d\omega_{m, \epsilon}(y)  \Big|^2 \, d\tilde{ E }_{ m , \epsilon} .
\end{aligned}
\end{align}
We use the embedding of $L_1$ in $L_2$ on a compact interval and Cauchy-Schwartz Inequality:
\begin{align}
\begin{aligned}
&\int_{\tilde{E}_{ m, \epsilon}} \Big| \dfrac{1}{| D( 0, c_r \epsilon ) | } \int_{ D(0, c_r \epsilon) }  \varphi \left( \dfrac{\zeta}{ c_r  \epsilon   } \right) D_y u \, d\omega_{ m, \epsilon} (y ) \Big|^2  \, d\tilde{E}_{ m,\epsilon}
\\
&\leq \int_{\tilde{ E }_{ m, \epsilon}}\Big( \dfrac{1}{ | D( 0, c_r \epsilon ) | } \int_{ D( 0 , c_r \epsilon ) }| \varphi \left( \dfrac{\zeta}{ c_r  \epsilon   } \right) D_{y} u |^2 \, d\omega_{ m , \epsilon } (y) \Big)  \, d \tilde{E}_{ m, \epsilon}
\\
&\leq \dfrac{ || \varphi ( \zeta/ c_r \epsilon )  ||_{L_2(D(0 , c_r \epsilon ))}^2 }{ | D( 0, c_r \epsilon ) |} \int_{\tilde{E}_{ m, \epsilon}} | D_y u|^2 \, d\tilde{E}_{ m,  \epsilon} .
\end{aligned}
\end{align}

Lastly we apply Proposition \ref{P:tildeEme_I} on (\ref{E:2.112}).
$\Box$

\begin{proposition} \label{P:Poincare_I}
For $u \in H^1({E}_{m,\epsilon})$, the averaging operator ${P}_{ m, \epsilon }$ satisfies a Poincar\'e-type inequality:
\begin{equation} 
\label{E:2.112}
\int_{ {E}_{m,\epsilon}} | u - {P}_{ m, \epsilon } u |^2 \, d {E}_{m,\epsilon} \leq c \epsilon^{2   } \int_{{E}_{m,\epsilon}} |\nabla u |^2 \, d {E}_{m,\epsilon} .
 \end{equation}
\end{proposition}

\textbf{Proof:}
Beginning with the fibered space, a calculation of the difference squared on each cross-section gives
\begin{align}
\begin{aligned}
| u - \tilde{P}_{m, \epsilon}  u |^2 
&=   \bigg| \int_{\varpi_{m, \epsilon} (y)  } \dfrac{ f_{y, \zeta}  (z, r, \theta) }{ r} D_\zeta u(y, \zeta)  \, d\omega_{m, \epsilon} (y) \bigg|^2
\\
&\leq 
c \bigg| \int_{ \varpi_{m, \epsilon} (y) } \dfrac{ D_\zeta u(y, \zeta)}{r}  \, d\omega_{ m, \epsilon } ( y ) \bigg|^2
\leq  c'   R_y  D_\zeta u(y, \zeta)
\end{aligned}
\end{align}
where $R_y$ is the operator of the form of Lemma \ref{L:Riesz} on $L_2(\omega_{m, \epsilon}(y))$ (in this case, it is the convolution with $1/r$).
From (\ref{E:R}) the 
norm of $R_y$ is bounded by $c \epsilon$.
Lastly we apply Proposition \ref{P:tildeEme_I}.
$\Box$

\subsection{Bounding the Norm on the Uniformly Fattened Binding}\label{SS:ext_binding_I}

Having established the required estimations for an averaging operator on each stratum, we now need to combine these different averaging operators into a global one.
To do so, here we establish several propositions regarding the trace on the interface $\Gamma_{ k, m, \epsilon}$ between $M_{k,S,\epsilon}$ and $E_{m,\epsilon}$.

\begin{definition} \label{D:Te}
The \textbf{trace }or \textbf{restriction operator} from $ {M}_{k, S, \epsilon}$ to $ {\Gamma}_{k, m, \epsilon}$ is denoted $T_{k,m,\epsilon}$.

The trace operator from $ {E}_{m,\epsilon}$ to $ {\Gamma}_{k, m, \epsilon}$ is denoted $T_{m, k, \epsilon}$.
\end{definition}

The standard embedding theorem claims that the trace space 
\\
$T_{k , m, \epsilon} H^1( {M}_{k, S, \epsilon})$ is isomorphic to $H^{1/2}(  {\Gamma}_{k,m,\epsilon})$.
\footnote{The trace space of $\Omega$ restricted to $\Gamma$ is given by the norm:
\begin{equation} \label{tracenorm}
 ||v ||_{TH^1(\Omega)} := \inf_{u \in H^1(\Omega): \, u|_{\Gamma} = v} || u||_{H^1(\Omega)}
 \end{equation}
}
However, these $\epsilon$-dependent spaces are in general not uniformly equivalent as metric spaces as is seen in \cite{MaPo}.

\begin{definition} \label{D:slobo}
Let $\Gamma$ be an $n$-dimensional domain.
Then $[f]_{\Gamma}$ denotes the following seminorm
\begin{equation} \label{slobo}
 [f]^2_{\Gamma} = \int_{\Gamma \times \Gamma} \dfrac{ | f(x) - f(y) |^2 }{ |x-y |^{1 + n}} \, dx \, dy .
 \end{equation}
The $H^{1/2}$ norm is given by $||u||_{H^{1/2}(\Gamma)}^2 = || u ||_{L_2(\Gamma)}^2 + [ u ]^2_{\Gamma}$.
\end{definition}

Let us estimate the trace on the fattened bindings.
First, we state a result that connects Proposition \ref{P:Eme_partition_I} to a trace estimation.

\begin{lemma} \label{L:special_lip}
Let $\Omega$ be a special Lipschitz domain and let $\varphi$ be the associated graph function with bounded Lipschitz norm $c_\Omega$.
Let $T_{\varphi}$ denote the operator from $L_2(\Omega)$ to $L_2(\R_+^n)$ (the half-space)
given by
\begin{equation}
T_{\varphi} u = u( x', x_n + \varphi( x') ) \qquad  x = (x', x_n) \in \R^n_+ .
\end{equation}
Then $T_\varphi$ is also a bounded linear operator from $H^1(\Omega)$ to $H^1(\R^n_+)$ whose norm depends only on the $\varphi$ and in particular $c_\Omega$.
\end{lemma}

\textbf{Proof:}
We begin with calculating the derivative (for $i < n$):
\begin{equation}
D_{x_i} T_{\varphi} u =   u_{x_i}(x', x_n + \varphi(x') ) +  u_{x_n}(x', x_n + \varphi(x')) D_{x_i} \varphi(x')
\end{equation} 
and (for $i = n$)
\begin{equation}
D_{x_n} T_{\varphi} u = u_{x_n} (x', x_n + \varphi(x')) .
\end{equation}
The Jacobian of the transformation $(x', x_n) \mapsto (x', x_n + \varphi(x'))$ also only depends on $\varphi$ and its derivatives.
Consequentially, the norm $T_{\varphi}: H^1(\Omega) \mapsto H^1(\mathbb{R}^n_+)$ has an upperbound that depends only on the maximum of $|\varphi|$ and $c_\Omega$. $\Box$



\begin{lemma} \label{L:trace_Eme_I}
Let $\{{E}_{m,\epsilon}\}$ be a family of fattened bindings ($\epsilon \in (0, \epsilon_0 ] $).
Let $u \in H^1( {E}_{m,\epsilon})$,
then one has
\begin{align}
\label{E:2.123}
\begin{aligned}
  \epsilon^{-1} || T_{m, k, \epsilon} & ( u -   P_{m, \epsilon }    u) ||^2_{L_2( {\Gamma}_{k,m,\epsilon })}
  + [ T_{ m, k, \epsilon}(u -  P_{m , \epsilon }  u) ]^2_{ {\Gamma}_{k, m,\epsilon }}
  \\
&  \leq c_m || u ||^2_{H^1( {E}_{m, \epsilon})} .
 \end{aligned}
 \end{align}
 \end{lemma}

 \textbf{Proof:}
We apply the partition of unity $\{\varphi_{i,\epsilon}\}$ as laid out in Proposition \ref{P:Eme_partition_I} and use Lemma \ref{L:special_lip} in the scaled domain.
We denote the homothetic scaling on $\R^3$: $\theta: x \to  x/\epsilon $ and $\Theta$ the induced operator on functions ($\Theta u = u(\theta)$).
Beginning with the left hand side of (\ref{E:2.123}), we have:
\begin{align}
\label{E:zero_ave}
\begin{aligned}
 \epsilon^{-1} || T_{m , k, \epsilon} &( u -    P_{m, \epsilon}   u) ||^2_{L_2( {\Gamma}_{k,m,\epsilon })}
 %
 + [ T_{ m, k,\epsilon}(u -   P_{m, \epsilon}   u) ]^2_{ {\Gamma}_{k, m, \epsilon }}
\\
& =
\epsilon^{ } || \Theta T_{m, k,\epsilon} ( u -    P_{m, \epsilon}   u) ||^2_{L_2(\theta( {\Gamma}_{k,m,\epsilon }))}
 \\
 &+ \epsilon^{ } [ \Theta T_{m,k, \epsilon}(u -   P_{m, \epsilon}   u)  ]^2_{\theta( {\Gamma}_{k, m, \epsilon })}
\\
& \leq
 \epsilon^{ }   \sum_i || \varphi_{i,\epsilon}(\theta) \Theta T_{m, k,\epsilon} ( u -   P_{m, \epsilon}   u)    ||^2_{L_2(\theta( {\Gamma}_{k, m, \epsilon }))}
\\
&+  [\varphi_{i,\epsilon}(\theta)  \Theta T_{m, k, \epsilon}(u -   P_{m, \epsilon}   u)     ]^2_{\theta( {\Gamma}_{k, m, \epsilon })}
.
\end{aligned}
\end{align}
Recalling Proposition \ref{P:Eme_partition_I}, we note $|\nabla \varphi_{i, \epsilon}(\theta)| \leq c_{  \varphi}$
and  $\mathsf{supp}(\varphi_{i, \epsilon}(\theta)) \leq c_0$ with both bounds uniform with respect to $\epsilon$.
Local finiteness of the partition holds as well (Proposition \ref{P:Eme_partition_I} (2)).
We identify $\mathsf{supp}(\varphi_{i, \epsilon})$ ($U_{i, \epsilon}$) with a local neighborhood of a special Lipschitz domain $\Omega_{i, \epsilon}$ with a graph norm bounded above by $c_M$ (see Proposition \ref{P:Eme_partition_I} (4)).

We then define $T_{\phi_{i, \epsilon}}: L_2(U_{i, \epsilon}) \mapsto L_2(\mathbb{R}^3)$ 
and denote the coordinate transformation from $U_{i, \epsilon}$ to $\mathbb{R}^3_+$ as $\chi_{i, \epsilon}$.
Subsequently, on each copy of $\R^3_+$, we invoke the Sobolev embedding theorem:

\begin{align}
\label{E:B.17}
\begin{aligned}
 ||  \Theta & T_{m, k,\epsilon} ( u
  -   P_{m, \epsilon}   u)    ||^2_{L_2(\mathsf{supp}(\varphi_{i, \epsilon}(\theta) ) \cap \theta( {\Gamma}_{k, m, \epsilon }))}
\\
&+ [ \Theta T_{m, k, \epsilon}(u -   P_{m, \epsilon}   u)     ]^2_{\mathsf{supp}(\varphi_{i, \epsilon}(\theta) ) \cap \theta( {\Gamma}_{k, m, \epsilon })}
\\
&\leq 
c \big(  ||  \Theta T_{\phi_{i, \epsilon}} T_{m, k,\epsilon} ( u -   P_{m, \epsilon}   u)    ||^2_{L_2(\mathsf{supp}(\varphi_{i, \epsilon}(\theta( \chi_{i, \epsilon}) ) ) \cap \theta(\chi_{i, \epsilon}( {\Gamma}_{k, m, \epsilon })))}
\\
&+  [   \Theta T_{\phi_{i, \epsilon}} T_{m, k, \epsilon}(u -   P_{m, \epsilon}   u)     ]^2_{\mathsf{supp}(\varphi_{i, \epsilon}(\theta(\chi_{i, \epsilon}))) \cap \theta( \chi_{i, \epsilon}({\Gamma}_{k, m, \epsilon }))} \big)
.
\end{aligned}
\end{align}
Denoting the upper bound of the norm of the embedding as $c_{em}$ (depending only on $c_M$, the upper bound on the Lipschitz norms of the boundary graphs), the right hand side of (\ref{E:B.17}) is bounded by:
\begin{equation}
c_{em}   ||  \Theta T_{\phi_{i, \epsilon}}( u -   P_{m, \epsilon}   u)    ||^2_{H^1(\mathsf{supp}(\varphi_{i, \epsilon}(\theta( \chi_{i, \epsilon}) ) }
.
\end{equation}
Thus (\ref{E:zero_ave}) is bounded by
\begin{equation}
\label{E:B.19}
\epsilon^{   } c' \sum_i ||\varphi_{i,\epsilon}(\theta)  \Theta (u -   P_{m, \epsilon}    u)  ||^2_{H^1(\mathsf{supp}(\varphi_{i,\epsilon}(\theta)))}
.
\end{equation}
After imputing all the constants associated with our partition of unity, (\ref{E:B.19}) is bounded by
\begin{equation}
 \epsilon^{  } c' c_U  ( 1 + c_{\nabla \varphi}) || \Theta    (u -     P_{m, \epsilon}    u)    ||^2_{H^1(\theta(  {E}_{m,\epsilon}) )}
 .
\end{equation}
Lastly, we scale the domain back to $\epsilon$ size to get the bound
$c || u ||^2_{H^1(  {E}_{m,\epsilon})} $.
$\Box$

With a norm estimate on the trace space of $ {E}_{m,\epsilon}$, we can now construct an extension operator from $ {\Gamma}_{k , m, \epsilon}$ to $ {M}_{k, S, \epsilon}$.

\begin{proposition} \label{P:ext_zero_average_I}
For $u \in H^1( {E}_{m,  \epsilon})$, the complement of the cross-sectional average $u -  P_{m, \epsilon}   u$ has an extension into ${M}_\epsilon$ denoted
$\mathcal{E}_{m,\epsilon} (u - P_{m, \epsilon}u)$  such that
\begin{equation} \label{E:u-p_me}
|| \mathcal{E}_{m,\epsilon} ( u -   P_{m, \epsilon}   u) ||^2_{H^1( {M}_\epsilon)} \leq c_m || u||^2_{H^1( {E}_{m,\epsilon})} .
\end{equation}
Furthermore, $\mathcal{E}_{m,\epsilon} (u -  P_{m, \epsilon} u)$ is supported within an $O(\epsilon^{ })$ neighborhood of $E_m$.
\end{proposition}

\textbf{Proof:}
This inequality follows from the partitioning ($U_{i, \epsilon}$, $\varphi_{i, \epsilon}$) and scaling $\theta$ seen in the proof of Lemma \ref{L:trace_Eme_I} and the Calder\'on-Stein Theorem (Theorem \ref{stein}).
In brief, let $\mathcal{E}_{m, 1}$ denote the extension operator in the sense of Theorem \ref{stein} on the $\theta(E_{m, \epsilon})$.
It follows
\begin{equation}
|| \Theta(u - P_{m, \epsilon} u) ||^2_{H^1(\theta(E_{m, \epsilon}))} 
\leq c \epsilon^{-1} || u - P_{m, \epsilon} u ||^2_{H^1(E_{m, \epsilon})}
.
\end{equation}

Let $\psi$ denote a smooth function which is identically $1$ on $\theta(E_{m, \epsilon})$ and $0$ outside distance $1$ from $\theta(E_{m, \epsilon})$.
Then we may define an extension of $u - P_{m, \epsilon}u$:
\begin{equation}
\mathcal{E}_{m, \epsilon}(u - P_{m, \epsilon}u ) = \Theta^{-1} \psi \mathcal{E}_{m, 1} \Theta ( u - P_{m, \epsilon}u)
.
\end{equation}
The inequality (\ref{E:u-p_me}) follows. $\Box$

\begin{corollary} \label{C:PmNTm_I}
For $u \in H^1( {E}_{m, \epsilon})$, one has:
\begin{equation}
\label{E:2.125}
 ||   P_{m, \epsilon }     u -     T_{k , m, \epsilon}  N_{k, \epsilon }  u ||_{L_2( {E}_{m,\epsilon})}^2  
 \leq c \epsilon^2 ||u||^2_{ H^1( {E}_{m, \epsilon})} .
 \end{equation}
\end{corollary}

\textbf{Proof:}
While $ T_{k, m, \epsilon} N_{k, \epsilon } u$ is a function on the interface $ {\Gamma}_{k, m,\epsilon}$, 
it is a constant on cross-sections $\partial \omega_{m, \epsilon}(x)$.
With an abuse of notation, we can set $N_{k, \epsilon} u(x \in E_m):=N_{k, \epsilon } u|_{\partial \omega_{m, \epsilon}(x)}$.
Beginning with an application of Proposition \ref{P:tildeEme_I},
we have
\begin{align}
		\label{E:2.126}
		\begin{aligned}
		||  \Phi_{E_{ m, \epsilon } }^{-1} & \tilde{P}_{m, \epsilon }   \Phi_{E_{ m, \epsilon } }   u -     T_{k , m, \epsilon}  N_{k, \epsilon }  u ||_{L_2( {E}_{m,\epsilon})}^2
		\\
		&\leq
		(1 + O(\epsilon) )
		||    \tilde{P}_{m, \epsilon }  \Phi_{E_{ m, \epsilon } }   u -     \Phi_{E_{ m, \epsilon } } T_{k , m, \epsilon}  N_{k, \epsilon }  u ||_{L_2( \tilde{E}_{m,\epsilon})}^2
		\\
		& =  (1 + O(\epsilon) ) \int_{E_m} \int_{\omega_{m,\epsilon} (y) } | \tilde{P}_{m, \epsilon }  \Phi_{E_{ m, \epsilon } }   u
		 \\
		 &-  \Phi_{E_{ m, \epsilon } } T_{k , m, \epsilon}  N_{k, \epsilon }  u|^2 \,   d\omega_{m, \epsilon}(y) d{E}_{m}.
		 \end{aligned}
 \end{align}
 Noting $\tilde{P}_{m, \epsilon} \Phi_{E_{m, \epsilon} } u$ can be extended to the boundary of $\tilde{E}_{m, \epsilon}$, (\ref{E:2.126}) is bounded by
\begin{align}
 \label{E:2.127}
 \begin{aligned}
\dfrac{ \max_{y \in E_m} |\omega_{m, \epsilon}(y)|}{2 \epsilon}  &\int_{\tilde{\Gamma}_{k, m, \epsilon } }  | \tilde{N}_{k , \epsilon} [ \tilde{P}_{m, \epsilon }  \Phi_{E_{ m, \epsilon } }   u
\\
 &-  \Phi_{E_{ m, \epsilon } } T_{k , m, \epsilon}    u ]  |^2 \,   d \tilde{\Gamma}_{k, m, \epsilon} .
\end{aligned}
\end{align}
Because the norm of $\tilde{N}_{k,\epsilon}$ is bounded independently of $\epsilon$, the above (\ref{E:2.127}) is bounded by
 \begin{equation}
 \label{E:2.128}
  c \epsilon  \int_{\tilde{\Gamma}_{k, m, \epsilon } }  | \tilde{P}_{m, \epsilon }  \Phi_{E_{ m, \epsilon } }   u
 -  \Phi_{E_{ m, \epsilon } } T_{k , m, \epsilon}   u|^2 \,   d \tilde{\Gamma}_{k, m, \epsilon} .
\end{equation} 
After applying the operator $\Phi_{E_{m, \epsilon} }^{-1}$, we have (\ref{E:2.128}) is equal to 
\begin{equation}
(1 + O(\epsilon) ) ||  {P}_{m, \epsilon }  \Phi_{E_{ m, \epsilon } }   u
 -   T_{k , m, \epsilon}   u ||^2_{L_2(\Gamma_{k, m , \epsilon})}
\end{equation}
This is the $L_2$ term in (\ref{E:2.123}), so we use Lemma \ref{L:trace_Eme_I}.
Consequentially, the desired bound for (\ref{E:2.125}) is achieved.
$\Box$

\begin{lemma} \label{L:Tkm_I}
For $u \in H^1( {M}_{k,S,\epsilon})$, one has:
\begin{equation}
  || T_{k, m, \epsilon}  u ||^2_{L_2( {\Gamma}_{k,m,\epsilon}) }  \leq c_k || u ||^2_{H^1( {M}_{k,   S,\epsilon}) } .
  \end{equation}
\end{lemma}

\textbf{Proof:}
The uniformly fattened page $M_{k, S, \epsilon}$ is a ``slab'' of width $2 \epsilon$.
We cover a neighborhood in $M_{k, S, \epsilon}$ of the interface $\Gamma_{k,m,\epsilon}$ with a partition of unity similar to Proposition \ref{P:Eme_partition_I} with some differences.
Let $\{U_{i, \epsilon} \}$ be collection of locally finite open covers of $\Gamma_{k, m, \epsilon}$ such that the maximum number of nontrivial intersections is bounded above by $n_U$ for all $\epsilon > 0$.
We also suppose the intersection $U_{i, \epsilon} \cap U_{i', \epsilon}$ contains a set of diameter larger than $c_1 \epsilon$.
The inner and outer diameters of each $U_{i, \epsilon}$ have lower and upper bounds of $c_2 \epsilon$ and $c_3 \epsilon$ respectively.

We consider cylindrical domains  $V_{i, \epsilon}$ in $M_{k, S, \epsilon}$ with $U_{i, \epsilon}$ as its base.
For some point $y \in U_{i, \epsilon}$, we denote the normal vector to $\Gamma_{k, m, \epsilon}$ at $y$ pointing into $M_{k, S, \epsilon}$ as $\nu_{k, m, \epsilon}(y)$.
For some constant $c_4$ (depending only on the geometry of $M_{k, S}$), the collection of sets $\{ V_{i, \epsilon} \}$ where
\begin{equation}
V_{i, \epsilon} := \{ x \in M_{k, S, \epsilon} : x = y + z \nu_{k, m, \epsilon}(y)  \qquad y \in U_{i, \epsilon},
\: z \in (0, c_4) \}
\end{equation}
has the finite intersection property as $\epsilon \to 0$. 
I.e. there is an $n_V$ such that at most $n_V$ sets $V_{i, \epsilon}$ (for a collection of $i$) have non-trivial intersection.

We equip $M_{k, S, \epsilon}$ with a local normal coordinate system $(y_1, y_2, z)$, where $y_2$ denotes the distance from the boundary $\Gamma_{k, m, \epsilon}$.
Considering the scaling $\theta: (y_1, y_2, z) \to ( y_1 / \epsilon, y_2, z/ \epsilon)$ (and induced operator $\Theta$).
Let $\{ \varphi_{i, \epsilon} \}$ be a smooth partition of unity subordinate to $\{ V_{i, \epsilon} \}$.
 Applying the scaling $\theta$ and the partition of unity, we have
\begin{equation}
\label{E:B.29}
   || T_{k , m, \epsilon}  u ||^2_{L_2( {\Gamma}_{k,m,\epsilon}) }
 \leq
 \sum_i  \epsilon^{2} ||  T_{k , m, \epsilon} \Theta  \varphi_{i,\epsilon}(\theta)   u  ||^2_{L_2(\theta( {\Gamma}_{k,m,\epsilon})) }
 .
 \end{equation}
 Under the scaling $\theta$, the support sets $ \mathsf{supp}( \varphi_{i, \epsilon}(\theta) )$ are contained in a ball of radius $c_1$ uniform with respect to $i$ and $\epsilon$ and contain a ball of radius $c_2$ also uniform with respect to $i$ and $\epsilon$.
 As before, the each of these domains is equivalent to a subset of a special Lipschitz domain whose graph function has Lipschitz norm bounded above by $M$ (also a uniform constant).
 The right hand side of (\ref{E:B.29}) is bounded by:
 \begin{equation}
  c_{em} \epsilon^{2} || \Theta  \varphi_{i,\epsilon}(\theta)    u  ||^2_{H^1(\theta( {M}_{k,S,\epsilon}))}
 \leq
c  ||    u  ||^2_{H^1(  {M}_{k, S,\epsilon} )}
.
\qquad \Box
\end{equation}

\begin{corollary} \label{C:NTkm_I}
For $u \in H^1({M}_\epsilon)$, one has:
\begin{equation}
||  T_{k,m, \epsilon} N_{k, \epsilon} u ||^2_{L_2( {E}_{m,\epsilon})} \leq c \epsilon^{ } || u||^2_{H^1(  {M}_{k, S,\epsilon})} .
\end{equation}
\end{corollary}

\textbf{Proof:}
It is analogous to the proof of Corollary \ref{C:PmNTm_I} using Lemma \ref{L:Tkm_I}.

\begin{theorem}  \label{T:Esmall_I}
For $u  \in H^1( {M}_\epsilon)$, the $L_2$ norm of $u$ on $ {E}_{m, \epsilon}$ is small:
\begin{equation}
 | |u | |^2_{L_2( {E}_{m, \epsilon})} \leq  c \epsilon^{ } | | u | |^2_{H^1( {M}_{\epsilon})} .
\end{equation}
\end{theorem}

\textbf{Proof:}
We use the triangle inequality:
\begin{align}
\begin{aligned}
 | | u  | &|_{L_2( {E}_{m,\epsilon})}
\leq
 || u -    P_{m, \epsilon }     u   ||_{L_2( {E}_{m,\epsilon})}
\\
&+  ||    P_{m, \epsilon }   u -     T_{k , m, \epsilon}  N_{k, \epsilon }  u ||_{L_2( {E}_{m,\epsilon})}
+
|| T_{k , m, \epsilon}  N_{k, \epsilon }  u ||_{L_2( {E}_{m,\epsilon})} 
.
\end{aligned}
\end{align}
With Proposition \ref{P:Poincare_I} and Corollaries \ref{C:PmNTm_I} and \ref{C:NTkm_I}, the theorem is proven. $\Box$

\begin{corollary} \label{C:H1Eme_I}
Assuming $u \in \mathcal{P}_{\Lambda}^\epsilon L_2(  {M}_\epsilon)$ for $\Lambda \leq c   \epsilon^{-1 + \delta }$, $\delta > 0$, and $\Lambda \notin \sigma(A_\epsilon)$, then the $H^1$-norm of $u$ on $ {E}_{m,\epsilon}$ is $o(1)$ with respect to $H^1$-norm on $M_\epsilon$.
\end{corollary}

\textbf{Proof:}
Due to the embedding of $\mathcal{P}_{\Lambda}^\epsilon L_2( {M}_\epsilon)$ into $L_2(  {M}_\epsilon)$, we can write
\begin{align}
\begin{aligned}
 |  | \nabla u  |   |^2_{L_2(  {E}_{m,\epsilon})}
& \leq \Lambda   | | u | |^2_{L_2(  {E}_{m,\epsilon})}
\\
& \leq c \Lambda \epsilon^{ } |  | u  | |_{H^1( {M}_\epsilon)}^2 = c \epsilon^\delta || u||^2_{H^1(  {M}_\epsilon)}
 \qquad \Box
\end{aligned}
\end{align}

\subsection{Averaging Operator $J_\epsilon$}\label{SS:exp_contr_I}

At last we can then define the averaging and extension operators in the sense of Definition \ref{D:Je}.

\begin{lemma} \label{L:d}
For any complex numbers $a$ and $b$ and for $d \in (0,1)$, one has:
\begin{equation}
\label{E:2.133}
(1 - d) |a|^2 + (1- d^{-1}) |b|^2 \leq | a + b |^2 \leq (1 + d) |a|^2 + (1 + d^{-1} ) |b|^2
.
\end{equation}
\end{lemma}

\textbf{Proof:}
Let us first assume both $a$ and $b$ are real. 
Because $(d^{1/2}  a \pm d^{-1/2}b)^2$ is non-negative,
\begin{equation}
-d a^2 - d^{-1} b^2
\leq 2 a b 
\leq d a^2 + d^{-1} b^2
.
\end{equation}
This completes the argument for the real case. 
The complex case follows from elementary arguments. 
$\Box$

\begin{Pro} \label{P:Je_I}
Let $M$ be an open book domain (Definition \ref{D:M}) and $M_\epsilon$ be the corresponding uniformly fattened domain  (Definition \ref{D:Me}).
Assume $\Lambda  \leq c \epsilon^{-1 + \delta}$ where $\delta > 0$ and $\Lambda \notin \sigma(A_\epsilon)$.
For some $\epsilon_0 > 0$, the family of linear operators $\{J_{\epsilon}  \}_{\epsilon \in (0, \epsilon_0] }$ that satisfies the conditions in Definition \ref{D:Je} for the open book structure $M$ is ($u \in \mathcal{P}_\Lambda^\epsilon L_2(M_\epsilon)$)
\begin{equation}
J_\epsilon  u := \begin{cases}
\sqrt{2\epsilon} \Psi_{M_k}^{-1}   N_{k, \epsilon} \big[ u + \sum_m \mathcal{E}_{m,\epsilon} (   P_{m , \epsilon }   u - u) \big] & M_{k, S, \epsilon } \mapsto M_{k }
\\
\sqrt{2\epsilon}  P_{m, \epsilon}   u
& E_{m, \epsilon} \mapsto E_m     
\end{cases}
\end{equation}
\end{Pro}

\textbf{Proof:}
First, we check whether $J_\epsilon u$ satisfies the boundary conditions on $\mathcal{G}^1$.
\begin{align}
\begin{aligned}
\lim_{x' \to x \in \partial M_{k, S, \epsilon} \cap \partial E_{m, \epsilon} } 
&N_{k, \epsilon}  \big[ u (x')+ \sum_m \mathcal{E}_{m,\epsilon} (  P_{m , \epsilon }    u - u)(x') \big]  
\\
&= N_{k, \epsilon}   P_{m , \epsilon }    u (x) = P_{m , \epsilon }   u (x)
.
\end{aligned}
\end{align}
Thus $J_\epsilon u$ is in $\mathcal{G}^1$.
Because each $\mathcal{E}_{m, \epsilon}( u - P_{m, \epsilon} u )$ is supported in a small $O(\epsilon)$ neighborhood around $E_m$, these extensions have disjoint supports.
Using Lemma \ref{L:d}, we break up the terms on $M_{k,S,\epsilon}$,
\begin{align}
\begin{aligned}
(1 - d) \big| \sqrt{2\epsilon}  & \Psi_{M_k}^{-1}  N_{k, \epsilon}   u   \big|^2\
\\
&+ 
(1 - d^{-1}) \sum_m \big| \sqrt{2\epsilon}  \Psi_{M_k}^{-1}  N_{k, \epsilon}     \mathcal{E}_{m,\epsilon} (  P_{m , \epsilon }   u - u)   \big|^2
\\
&\leq
\big| \sqrt{2\epsilon}   \Psi_{M_k}^{-1} N_{k, \epsilon}  \big[ u + \sum_m \mathcal{E}_{m,\epsilon} (  P_{m , \epsilon }   u - u) \big] \big|^2
\\
&\leq
(1 + d) \big| \sqrt{2\epsilon}   \Psi_{M_k}^{-1}  N_{k, \epsilon}   u   \big|^2\
\\
&+ 
(1 + d^{-1}) \sum_m \big| \sqrt{2\epsilon}  \Psi_{M_k}^{-1}  N_{k, \epsilon}     \mathcal{E}_{m,\epsilon} (  P_{m , \epsilon }   u - u)   \big|^2 .
\end{aligned}
\end{align}
To demonstrate the $L_2$ near isometry property, we first assume that
 $||J_\epsilon u ||^2_{L_2(M_{k} )} \geq ||u||^2_{L_2(M_{k, S, \epsilon})}$. 
The demonstration of the 
other case $||J_\epsilon u ||^2_{L_2(M_{k} )} \leq ||u||^2_{L_2(M_{k, S, \epsilon})}$ 
requires only minor modification.
We calculate the upper and lower bound on the norm difference:
\begin{align}
\label{E:2.138}
\begin{aligned}
\sum_k  
(1 &- d)  || \sqrt{2\epsilon}   \Psi_{M_k}^{-1}  N_{k, \epsilon}   u   ||^2_{L_2(M_{k })} 
\\
&+ (1 - d^{-1})
 \sum_{k, m} || \sqrt{2\epsilon}  \Psi_{M_k}^{-1}  N_{k, \epsilon}     \mathcal{E}_{m,\epsilon} (  P_{m , \epsilon }   u - u)   ||^2_{L_2(M_{k })} 
\\
&- \sum_k ||u ||_{L_2(M_{k, S, \epsilon})}^2   - ||u||^2_{L_2(E_{m, \epsilon})}
\\
&\leq
   ||J_\epsilon u ||_{L_2(M)}^2 - || u ||_{L_2(M_\epsilon)}^2  
\\
&\leq
\sum_k  
(1 + d) || \sqrt{2\epsilon}   \Psi_{M_k}^{-1}  N_{k, \epsilon}   u   ||^2_{L_2(M_{k, S, \epsilon})} 
\\
&+ (1 + d^{-1})
 \sum_{k, m} || \sqrt{2\epsilon}  \Psi_{M_k}^{-1}  N_{k, \epsilon}    \mathcal{E}_{m,\epsilon} (  P_{m , \epsilon }   u - u)   ||^2_{L_2(M_{k, S, \epsilon})} 
\\
&- \sum_k ||u ||_{L_2(M_{k, S, \epsilon})}^2   - ||u||^2_{L_2(E_{m, \epsilon})} .
\end{aligned}
\end{align}
Since we only require demonstrating that $||J_\epsilon u||_{H^1(M)}$ is bounded above (\ref{E:J2}), we begin with assuming  $||J_\epsilon u||_{H^1(M)} \geq || u ||_{H^1(M_\epsilon)}$ and write:
\begin{align}
\label{E:2.139}
\begin{aligned}
   ||J_\epsilon u &||_{H^1(M)}^2 
   - || u ||_{H^1(M_\epsilon)}^2  
\\
&\leq
\sum_k  
(1 + d) || \sqrt{2\epsilon}   \Psi_{M_k}^{-1}  N_{k, \epsilon}   u   ||^2_{H^1(M_{k, S, \epsilon})} 
\\
&+ (1 + d^{-1})
 \sum_{k, m} || \sqrt{2\epsilon}  \Psi_{M_k}^{-1}  N_{k, \epsilon}    \mathcal{E}_{m,\epsilon} (  P_{m , \epsilon }   u - u)   ||^2_{H^1(M_{k, S, \epsilon})} 
\\
&- \sum_k ||u ||_{H^1(M_{k, S, \epsilon})}^2   - ||u||^2_{H^1(E_{m, \epsilon})} 
.
\end{aligned}
\end{align}

Having established these two inequalities (\ref{E:2.138}) and (\ref{E:2.139}),
we collect terms in these inequalities and apply various propositions established in this chapter to demonstrate which terms are negligible (are $o(1)$ in an appropriate norm) and which terms are nearly an isometry (are $1 + o(1)$ in an appropriate norm).

By Proposition \ref{P:MkstoMk_I}, we have
\begin{align}
\begin{aligned}
\big| \int_{M_{k }} & | \sqrt{2\epsilon}  \Psi_{M_k}^{-1}  N_{k, \epsilon}  u |^2 \, dM_k   
- \int_{M_{k, S, \epsilon}} |u |^2 \, dM_{ \epsilon} \big|
\\
&\leq
\big| (1 + O(\epsilon))\int_{M_{k, S}} | \sqrt{2\epsilon}     N_{k, \epsilon}  u |^2 \, dM_k   - \int_{M_{k, S, \epsilon}} |u |^2 \, dM_{ \epsilon} \big|
\\
&\leq
c \epsilon^{ } || u ||^2_{H^1( M_\epsilon) }
\end{aligned}
\end{align}
where the last inequality results from Proposition \ref{P:Nu_I}.
We note the energy bound only needs to be demonstrated from above, so we see
\begin{align}
\begin{aligned}
 \int_{M_{k }} |  \nabla_{M_k} \sqrt{2\epsilon} & \Psi_{M_k}^{-1}  N_{k, \epsilon}  u |^2 \, dM_k   - \int_{M_{k, S, \epsilon}} | \nabla u |^2 \, dM_{ \epsilon} 
 \\
&\leq c \epsilon || u ||^2_{H^1( M_\epsilon) }
\end{aligned}
\end{align}
which follows from Propositions \ref{P:MkstoMk_I} and \ref{P:QNu_I}.

 This leaves the extensions from the fattened bindings into the page ($\mathcal{E}_{m,\epsilon}( u - P_{m, \epsilon} u)$) and the norm of the binding unaccounted for in (\ref{E:2.138}) and (\ref{E:2.139}).
 We estimate the $H^1$-norm of the extensions.
Using Propositions \ref{P:MkstoMk_I}, \ref{P:Nu_I}, and \ref{P:QNu_I}, and the disjoint supports of $E_{m, \epsilon} (u - P_{m, \epsilon} u )$:
\begin{equation} \label{E:J_i}
\begin{split}
 \sum_m || \sqrt{2\epsilon}   \Psi_{M_k}^{-1}   N_{k, \epsilon}  \mathcal{E}_{m,\epsilon} (   P_{ m, \epsilon}   u - u) \big  ||_{H^1(M_{k, S} ) }^2 + \sum_m ||u ||^2_{H^1(E_{m, \epsilon})}
\\
\leq
(1 + O(\epsilon) )   \sum_m || \mathcal{E}_{m,\epsilon} (   P_{m, \epsilon}   u - u) \big  ||_{H^1(M_{k, S, \epsilon } ) }^2
+ \sum_m ||u ||^2_{H^1(E_{m, \epsilon})}
.
\end{split}
\end{equation}
By Proposition \ref{P:ext_zero_average_I}, this is bounded by
\begin{equation}
(1 + O(\epsilon) ) c \sum_m || u ||^2_{H^1 (E_{ m, \epsilon } ) } .
\end{equation}
Because $u \in \mathcal{P}_{\Lambda}^\epsilon L_2(M_\epsilon)$ and Corollary \ref{C:H1Eme_I}, we arrive to the following upper bound on the norm of (\ref{E:J_i}):
\begin{equation}
c \epsilon^{\delta}  || u ||^2_{H^1(M_\epsilon) } .
\end{equation}

Hence by setting $d = \epsilon^{\delta/2}$, we conclude that $J_\epsilon u|_{M_{k}}$ is close in $L_2$ to $u$ and $J_\epsilon u|_{M_{k}}$ does not exceed the energy on $M_{ \epsilon}$ by more than an $o(1)$ factor.

Thus $J_\epsilon$ is an averaging operator as required in Theorem \ref{T:main} completing the proof of Proposition \ref{P:Je_I} and consequentially Theorem \ref{T:main}. $\Box$

\section{Discussion and Conclusions}\label{S:remarks}

These results lay a foundation for exploring the properties of these fattened domains when the underlying space is greater than one dimension.
We recognize that many of the works of the fattened graph domains could be extended to more general fattened domains in $\mathbb{R}^3$.
Having established a methodology of proving spectral convergence in the simplest case,
we can continue to extend this work to consider non-uniformly fattened domains, resolvent convergence, scattering problems, and so on.
Such results would be motivated by the modeling of micro-electronic or photonic systems.

Preliminary work by the author shows the proof can be extended to allow non-uniformly fattened structures.
Let us consider varying the width of the fattened pages by some 
continuously differentiable function $w(x)$,
(i.e. fattening by balls of radius $\epsilon w(x)$).
This changes the limit operator $A$ to be the weighted Laplace-Beltrami operator.
A more interesting result is changing the speed of shrinkage of the strata relative
to one another, namely by fattening the binding by balls of radius $\epsilon^\beta$.
These results are non-trivial, result in phase transitions,
 and will be presented elsewhere.
Further preliminary results explore the problem of fattened polyhedral domains and show a similar phase transition based on a capacity heuristic.

This still leaves many avenues of research that could show non-trivial physical phenomena.
One should eventually consider fattened domains with equipped with Schr\"{o}dinger operators with unbounded potentials, 
periodic fattened domains,
and tangential contact between $2$-strata.


\section{Acknowledgments}\label{acknow}

The work of the author was partially supported by the NSF DMS-1517938 Grant.

\bibliography{uniform_fattened}

\bibliographystyle{plain}



\end{document}